\documentclass{m2an}
\usepackage{a4wide,enumerate,xcolor}
\usepackage{amsmath,comment}
\usepackage[utf8]{inputenc}
\usepackage[T1]{fontenc}
\allowdisplaybreaks
\usepackage{enumitem}
\usepackage{bm} 
\usepackage{graphicx}
\usepackage[normalem]{ulem}
\usepackage{cancel}
\setlist[enumerate]{label={\rm(\arabic*)}, leftmargin=28pt, itemsep=3pt}
\setlist[itemize]{label={$\bullet$}, leftmargin=12pt, itemsep=3pt}
\newcommand{\phm}{{\phantom{-}}}
\let\pa\partial
\let\na\nabla
\let\eps\varepsilon
\newcommand{\N}{{\mathbb N}}
\newcommand{\R}{{\mathbb R}}
\newcommand{\diver}{\operatorname{div}}
\newcommand{\red}{\textcolor{red}}

\def\softd{{\leavevmode\setbox1=\hbox{d}%
\hbox to 1.05\wd1{d\kern-0.4ex{\char039}\hss}}}

\newcommand{\pvec}[1]{\vec{#1}\mkern2mu\vphantom{#1}}

\newtheorem{theorem}{Theorem}
\newtheorem{lemma}[theorem]{Lemma}

\newtheorem{remark}[theorem]{Remark}

\begin{document}
\title{Error estimates for a structure-preserving finite-element scheme of quasi-incompressible Navier--Stokes--Maxwell--Stefan systems}
\thanks{The first author was supported by the Gutenberg Research College and the Deutsche Forschungsgemeinschaft (DFG), grants 233630050 (TRR 146) and 441153493 (SPP 2256). The second author acknowledges partial support from the Austrian Science Fund (FWF), grant 10.55776/PAT2687825, and from the Austrian Federal Ministry for Women, Science and Research and implemented by \"OAD, project MULT09/2025. The research of the third author was supported by the DFG, grants 233630050 (TRR 146) and 525800857 (SPP 2410). She gratefully acknowledges support of the Gutenberg Research College of the University of Mainz and the Mainz Institute of Multiscale Modeling. This work has received funding from the European Research Council (ERC) under the European Union's Horizon 2020 research and innovation programme, ERC Advanced Grant NEUROMORPH, no.~101018153. For open access purposes, the authors have applied a CC BY public copyright license to any author-accepted manuscript version arising from this submission.} 

\author{Aaron Brunk}\address{Institute of Mathematics, Johannes Gutenberg-University Mainz, Staudinger Weg 9, 55288 Mainz, Germany,
abrunk@uni-mainz.de  \& lukacova@uni-mainz.de}

\author{Ansgar Jüngel}\address{Institute of Analysis and Scientific Computing, TU Wien, Wiedner Hauptstra\ss e 8--10, 1040 Wien, Austria,
\email{juengel@tuwien.ac.at} }

\author{M{\'a}ria Luk{\'a}{\v{c}}ov{\'a}-Medvi{\v{d}}ov{\'a}}\sameaddress{1}

\date{\today}

\begin{abstract}
A conforming finite-element scheme for quasi-incompressible Navier--Stokes--Maxwell--Stefan systems in a bounded domain with periodic
boundary conditions is presented. The system consists of the Navier--Stokes equations, together with a quasi-incompressibility constraint, coupled with the cross-diffusion Maxwell--Stefan equations.
The numerical scheme preserves the partial masses and the quasi-incompressibility constraint and dissipates the discrete energy. First-order error estimates in space and time are shown for a diffusion-augmented mobility matrix by using the relative energy method and suitable projection error estimates. Numerical experiments in two space dimensions illustrate the convergence of the scheme and indicate second-order convergence in space in the $L^2$ norm.
\end{abstract}

 

\maketitle


\section{Introduction}

The quasi-incompressible Navier--Stokes--Maxwell--Stefan equations provide a continu\-um-scale framework for modeling the coupled transport of mass and momentum in fluid mixtures. Assuming that the mixture is volume-additive, the fluid velocity satisfies an incompressibility constraint, but it is generally not divergence-free. The Maxwell--Stefan relations describe the multicomponent mass diffusion driven by concentration and pressure and capture cross-diffusion effects beyond the standard Fick law. The mass diffusion fluxes interact with the momentum balance, enabling the system to represent complex processes such as reactive flows, multiphase mixing, and transport in non-ideal mixtures. 

The existence of weak solutions to this system was proved in \cite{Dru21}, and a structure-preserving finite-element scheme was proposed in \cite{BJL25} without providing a convergence analysis. In this paper, we show the convergence of the scheme and, additionally, prove first-order error estimates for the particle densities, chemical potentials, fluid velocity, and pressure.

The dynamics of the particle densities $\rho_i$ and fluid velocity $\bm{u}$ is governed by the mass and momentum balance equations
\begin{align}
  & \pa_t\rho_i + \diver(\rho_i\bm{u} + \bm{J}_i) = 0,
  \label{1.mass} \\
  & \pa_t(\rho\bm{u}) + \diver\big(\rho\bm{u}\otimes\bm{u}
  - \mathbb{S}(\vec\rho,\na\bm{u}) + p\mathbb{I}\big) = 0
  \quad\mbox{in }\Omega,\ t>0, \label{1.mom} \\
  & \rho_i(\cdot,0) = \rho_{i}^0, \quad 
  (\rho\bm{u})(\cdot,0) = \rho^0\bm{u}^0\quad\mbox{in }\Omega,\ 
  i=1,\ldots,N, \label{1.ic}
\end{align}
where $\Omega$ is the $d$-dimensional torus ($d\ge 1$), $\rho=\sum_{i=1}^N\rho_i$ is the total density, and $p$ is the pressure. We denote the $N$-dimensional vector $(\rho_1,\ldots,\rho_N)$ by $\vec\rho$ and similar for other variables; $d$-dimensional variables like $\bm{u}=(\bm{u}_1,\ldots,\bm{u}_d)$ are written in bold font. We choose the stress tensor
\begin{align*}
  \mathbb{S}(\vec\rho,\na\bm{u}) = \nu(\na\bm{u}+\na\bm{u}^T)
  + \lambda(\diver\bm{u})\mathbb{I},
\end{align*}
where $\nu=\nu(\vec\rho)>0$ and $\lambda=\lambda(\vec\rho)\ge -2\nu/d$ are the density-dependent viscosities and $\mathbb{I}$ is the unit matrix in $\R^{d\times d}$. By Korn's inequality, the stress tensor satisfies the coercivity condition $\int_\Omega\mathbb{S}(\vec\rho,\na\bm{u}):\na\bm{u}\,dx\ge c\int_\Omega|\na\bm{u}|^2dx$ for some $c>0$, where ``:'' is the Frobenius matrix product and $|\cdot|$ the Euclidean norm. 

The partial diffusion flux $\bm{J}_i$ is given by the linear combination of the gradients of the chemical potentials $\mu_j$,
\begin{align*}
  \bm{J}_i = -\sum_{j=1}^N M_{ij}(\vec\rho)\na\mu_j, \quad
  i=1,\ldots,N,
\end{align*}
where the mobility matrix $M=(M_{ij})\in\R^{N\times N}$ is symmetric and positive semi-definite. As the net diffusion flux $\sum_{i=1}^N\bm{J}_i=0$ vanishes by assumption, the mobility matrix has a nontrivial kernel. The chemical potential $\mu_i$ is defined by
\begin{align}\label{1.chem}
  \mu_i = \frac{\pa f}{\pa\rho_i}(\vec\rho) + V_ip, \quad i=1,\ldots,N,
\end{align}
where $f$ is the internal energy density and $V_i>0$ is the specific volume of the $i$th species of the mixture. An example for ideal gas mixtures is $f(\vec\rho)=\sum_{i=1}^N\rho_i\log(\rho_i/\rho)$. We assume only that the function $f$ is homogeneous of degree one. 

We assume that the mixture is volume-additive, i.e., the incompressibility constraint
\begin{align}\label{1.incom}
  \sum_{i=1}^N V_i\rho_i = 1 \quad\mbox{in }\Omega
\end{align}
is satisfied. Then definition \eqref{1.chem} is consistent with the Gibbs--Duhem relation
\begin{align*}
  -f(\vec\rho) + \sum_{i=1}^N\rho_i\mu_i
  = -\sum_{i=1}^N\rho_i\frac{\pa f}{\pa\rho_i}
  + \sum_{i=1}^N\rho_i\mu_i = \sum_{i=1}^N \rho_i V_i p = p,
\end{align*} 
where the first identity follows from the positive homogeneity of $f$ and the last identity is a consequence of \eqref{1.incom}. The incompressibility constraint shows that the velocity field is generally not divergence-free. Indeed, multiplying \eqref{1.mass} by $V_i$ and summing over $i=1,\ldots,N$ leads to
\begin{align*}
  0 = \pa_t\bigg(\sum_{i=1}^N V_i\rho_i\bigg) 
  + \diver\bigg(\sum_{i=1}^N V_i\rho_i\bm{u}\bigg) 
  + \diver\bigg(\sum_{i=1}^N V_i\bm{J}_i\bigg)
  = \diver\bm{u} + \diver\bigg(\sum_{i=1}^N V_i\bm{J}_i\bigg).
\end{align*}
In particular, a net local change of volume is possible because of mixing. If all specific volumes are the same, we recover the incompressible Navier--Stokes--Maxwell--Stefan equations analyzed in \cite{ChJu15}. 

We summarize the reformulated model equations:
\begin{align}
  & \pa_t\rho_i + \diver(\rho_i\bm{u})
  - \diver\bigg(\sum_{j=1}^N M_{ij}(\vec\rho)\na\mu_j\bigg) = 0, 
  \label{1.mass2} \\
  & \pa_t(\rho\bm{u}) + \diver\bigg\{\bm{u}\otimes\bigg(\rho\bm{u} 
  - \sum_{i,j=1}^N M_{ij}(\vec\rho)\na\mu_j\bigg)
  - \mathbb{S}(\vec\rho,\na\bm{u}) + p\mathbb{I}\bigg\} = 0, \label{1.mom2} \\
  & \mu_i = \frac{\pa f}{\pa\rho_i} + V_ip, \quad i=1,\ldots,N,
  \label{1.mu} \\
  & \diver\bm{u} = \diver\bigg(\sum_{i,j=1}^N V_iM_{ij}(\vec\rho)
  \na\mu_j\bigg)\quad\mbox{in }\Omega,\ t>0, \label{1.div}
\end{align}
together with the initial conditions \eqref{1.ic}. Observe that the sum in \eqref{1.mom2} vanishes; it is explicitly written since we choose below a diffusion-augmented mobility matrix yielding a non-zero discrete total diffusion flux; see \eqref{num.mom}. 

A structure-preserving finite-element method for \eqref{1.mass2}--\eqref{1.div} was proposed in \cite{BJL25}. This scheme preserves the partial masses and the constraint \eqref{1.incom}, conserves the total mass, and satisfies a discrete inequality version of the energy equality
\begin{align*}
  \frac{d}{dt}\int_\Omega E(\vec\rho,\bm{u})dx
  + \int_\Omega\bigg(\mathbb{S}(\vec\rho,\na\bm{u}):\na\bm{u}
  + \sum_{i,j=1}^N M_{ij}(\vec\rho)\na\mu_i\cdot\na\mu_j\bigg)dx = 0,
\end{align*}
where the total energy density $E(\vec\rho,\bm{u}) = \frac12\rho|\bm{u}|^2 + f(\vec\rho)$ is the sum of the kinetic and internal energy densities. The convergence of the scheme has remained an open problem. The main result of the paper is contained in Theorem \ref{thm.error}: For a given regular solution $(\vec\rho,\vec\mu,\bm{u},p)$ to \eqref{1.ic}, \eqref{1.mass2}--\eqref{1.div}, we prove convergence rates for the differences $\vec\rho_h-\vec\rho$, $\vec\mu_h-\vec\mu$, $\bm{u}_h-\bm{u}$, and $p_h-p$ in suitable norms. 

To prove this result, we introduce the augmented mobility matrix $M_\eps:=M+\eps\mathbb{I}$ for some $\eps>0$, which ensures the positive definiteness of $M_\eps$. As a consequence, the total diffusion flux does not vanish anymore. We include this term in the convection term to ensure the cancellation of some higher-order terms in the error estimates; see \eqref{eq:discrelC} below. For the proof of the error estimates, we split them into the evolution error  and the projection error. The latter one is shown by standard interpolation and projection estimates. The evolution error estimate is based on a nonlinear stability result, which is proved by the relative energy method, comparing the discrete solution with the projected solution. The key steps are the residual error estimates, which require some careful computations. 

Let us discuss briefly the state of the art. The existence of weak solutions to \eqref{1.ic}, \eqref{1.mass2}--\eqref{1.div} was established in \cite{BoDr21,Dru21}. Structure-preserving numerical approximations for transient Maxwell--Stefan equations without fluid coupling can be found, e.g., in \cite{CEM24,HLTW21,JuZu23,McBo14}. In \cite{AFMV25,VFM22}, numerical approximations of coupled Maxwell--Stefan fluid models are analyzed using finite-element methods, where the steady-state Maxwell--Stefan equations are coupled with a compressible Stokes flow. Similar steady-state models and discretizations in the quasi-incompressible case are studied in \cite{BaFa25}, which uses an explicit force formulation rather than a flux-based one. To include additional effects, like in \cite{FMV23}, the explicit force formulation may be more suitable. For a related but different model, the quasi-incompressible Navier--Stokes--Cahn--Hilliard system, energy dissipative finite-element discretizations were presented in \cite{BEHL23,GiPr15,GLL14}.

The paper is organized as follows. The numerical scheme and the main results are made precise in Section \ref{sec.main}. In Section \ref{sec.prop}, we recall the structure-preserving properties of the scheme. Section \ref{sec.est1} is devoted to the technical proof of the error estimates, and we discuss the case $\eps=0$ in Section \ref{sec.eps}. Finally, some numerical experiments are presented in Section \ref{sec.exp} confirming the theoretical convergence results.


\section{Numerical scheme and main results}\label{sec.main}

We formulate the numerical scheme, state our main results, and recall the structure-preserving results of \cite{BJL25}. First, we introduce some notation. We already mentioned that $\vec\rho=(\rho_1,\ldots,\rho_N)$ denotes an $N$-dimensional vector, $\bm{u}=(\bm{u}_1,\ldots,\bm{u}_d)$ is a $d$-dimensional vector, and similar for other variables. The bracket $\langle\cdot,\cdot\rangle$ denotes the duality product on $H^{-k}(\Omega)\times H^k(\Omega)$ for $k\in\N$ or the inner product of $L^2(\Omega)$, depending on the context. We write $\|\cdot\|_{k,p}$ for the norm of $W^{k,p}(\Omega)$ with $k\ge 0$, $1\le p\le\infty$, $\|\cdot\|_{0,p}$ for the norm of $L^p(\Omega)$, and $|\cdot|_{1,p}$ is the seminorm in $W^{1,p}(\Omega)$, i.e.\ $|v|_{1,p}=\|\na v\|_{0,p}$ for $v\in W^{1,p}(\Omega)$. 

\subsection{Numerical scheme}

We regularize the equations by replacing $M_{ij}(\vec\rho)$ by $M_{\eps,ij}(\vec\rho) := M_{ij}(\vec\rho) + \eps\delta_{ij}$ for $\eps\ge 0$. The advantage is that the regularized mobility matrix becomes positive definite. However, the total mass flux does not vanish anymore, and we need to add the expression $\eps\sum_{i=1}^N\na\mu_i$ to equation \eqref{1.mom2} to ensure the energy balance. This leads to the following regularized equations:
\begin{align}
  & \pa_t\rho_i + \diver(\rho_i\bm{u})
  - \diver\bigg(\sum_{j=1}^N 
  M_{\eps,ij}(\vec\rho)\na\mu_j\bigg) = 0, \label{2.mass} \\
  & \pa_t(\rho\bm{u}) + \diver\bigg\{\bm{u}\otimes\bigg(\rho\bm{u} 
  - \eps\sum_{j=1}^N \na\mu_j\bigg)
  - \mathbb{S}(\vec\rho,\na\bm{u}) + p\mathbb{I}\bigg\} = 0, 
  \label{2.mom} \\
  & \mu_i = \frac{\pa f}{\pa\rho_i} + V_ip, \quad i=1,\ldots,N, 
  \label{2.mu} \\
  & \diver\bm{u} = \diver\bigg(\sum_{i,j=1}^N V_iM_{\eps,ij}(\vec\rho)
  \na\mu_j\bigg). \label{2.div}
\end{align}
The variational formulation of \eqref{2.mass}, \eqref{2.mu}, and \eqref{2.div} is as usual, while the weak formulation of \eqref{2.mom} becomes, for smooth test functions $\bm{v}$ (see \cite[Lemma 4]{BJL25}),
\begin{align*}
  & \frac12\langle \bm{u}\pa_t\rho,\bm{v}\rangle
  + \langle\rho\pa_t\bm{u},\bm{v}\rangle
  + \bm{B}_{\rm skw}(\rho\bm{u};\bm{u},\bm{v}) 
  + \bm{B}^\eps_{\rm skw}(\bm{J};\bm{u},\bm{v}) \\
  &\phantom{xx}- \langle p,\diver\bm{v}\rangle
  + \langle\mathbb{S}(\vec\rho,\na\bm{u}),\na\bm{v}\rangle 
  + \sum_{i=1}^N\langle\rho_i\na(\mu_i-V_ip),\bm{v}\big\rangle = 0,
\end{align*}
where the skew-symmetric bilinear forms are defined by
\begin{align*}
  (\bm{v},\bm{w})\mapsto \bm{B}_{\rm skw}(\rho\bm{u};\bm{v},\bm{w})
  &= \frac12\sum_{i=1}^d\bigg\langle
  \rho \bm{u}_i,\frac{\pa\bm{v}}{\pa x_i}\cdot\bm{w}
  - \frac{\pa\bm{w}}{\pa x_i}\cdot\bm{v}\bigg\rangle, \\
  (\bm{v},\bm{w})\mapsto \bm{B}^\eps_{\rm skw}(\bm{J};\bm{v},\bm{w})
  &= \frac12\sum_{i=1}^d\bigg\langle
  J_i,\frac{\pa\bm{v}}{\pa x_i}\cdot\bm{w}
  - \frac{\pa\bm{w}}{\pa x_i}\cdot\bm{v}\bigg\rangle,
\end{align*}
where $\bm{J}=(J_1,\ldots,J_d):=-\eps\sum_{i=1}^N\na\mu_i$ is the total mass flux. 

For the numerical discretization of the variational formulation of \eqref{2.mass}--\eqref{2.div}, we define the following spaces on the partition $[t^{k-1},t^k]$ of the time interval $[0,T]$, where $t^k=k\tau$, $k=1,\ldots,n$, and $\tau=T/n$:
\begin{align*}
  \Pi_\tau^0(\mathcal{U}) 
  &= \mbox{space of (discontinuous) piecewise constant functions on }
  \{t^0,\ldots,t^n\} \\
  &\phantom{xx} \mbox{with values in }\mathcal{U}, \\
  \Pi_\tau^1(\mathcal{U})  
  &= \mbox{space of  continuous piecewise linear functions on }
  \{t^0,\ldots,t^n\} \\
  &\phantom{xx}\mbox{ with values in }\mathcal{U}.
\end{align*}
The spatial discretization $\mathcal{T}_h$ is a geometrically conforming partition of $\Omega$ into simplices, which can be extended periodically to periodic extensions of $\Omega$, where $h>0$ is the maximal diameter of the simplices. In particular, the boundary of $\Omega$ coincides with the boundary of the triangulation. We define the following spaces:
\begin{align*}
  & \mathcal{V}_h = \big\{v\in H^1(\Omega)\cap C^0(\overline\Omega):
  v|_K\in P_1(K)\mbox{ for all }K\in\mathcal{T}_h\big\}, \quad
  \mathcal{X}_h = \mathcal{V}_h^d, \\
  & \mathcal{Q}_h = \{v\in\mathcal{V}_h:\langle v,1\rangle=0\}, \quad
  \mathcal{V}_{h,+} = \{v\in\mathcal{V}_h:v(x)>0
  \mbox{ for }x\in\Omega\},
\end{align*}
where we have set, slightly abusing the notation, $\langle u,v\rangle = \sum_{K\in\mathcal{T}_h}\int_K u|_K v|_K\,dx$, and $P_1(K)$ is the space of affine functions on $K$. The space $\mathcal{X}_h$ in \cite{BJL25} was chosen as the space of (continuous) piecewise quadratic functions. Here, we choose the space of piecewise linear functions to simplify the computations of the residual error estimates. Although the pair $\mathcal{X}_h\times\mathcal{Q}_h$ is not inf-sup stable, this does not give any problems, since we show below that the pressure error is controlled via the chemical potentials and not via an inf-sup condition.

Let the initial data $(\pvec{\rho}_h^0,\bm{u}_h^0)\in\mathcal{V}_{h,+}^N\times\mathcal{X}_h$ be given. We wish to find $(\pvec{\rho}_h^k,\bm{u}_h^k)\in\Pi_\tau^1(\mathcal{V}_{h,+}^N\times\mathcal{X}_h)$ and $(\vec\mu_h^k,p_h^k)\in\Pi_\tau^0(\mathcal{V}_h^N\times\mathcal{Q}_h)$ recursively. For this, we introduce for $\bm{u}_h$, $\bm{v}_h\in\mathcal{X}_h$ the discrete skew-symmetric bilinear forms
\begin{align*}
  \bm{b}_{\rm skw}(\rho_h\bm{u}_h;\bm{u}_h,\bm{v}_h)
  &= \frac12\sum_{i=1}^d\bigg\langle\rho_h\bm{u}_{h,i},
  \frac{\pa\bm{u}_h}{\pa x_i}\cdot\bm{v}_h
  - \frac{\pa\bm{v}_h}{\pa x_i}\cdot\bm{u}_h\bigg\rangle, \\
  \bm{b}^\eps_{\rm skw}(\bm{J}_h;\bm{u}_h,\bm{v}_h)
  &= \frac12\sum_{i=1}^d\bigg\langle J_{h,i},
  \frac{\pa\bm{u}_h}{\pa x_i}\cdot\bm{v}_h
  - \frac{\pa\bm{v}_h}{\pa x_i}\cdot\bm{u}_h\bigg\rangle,
\end{align*}
recalling the total mass flux is given by
\begin{align*}
  \bm{J}_h = (J_{h,1},\ldots,J_{h,d}) 
  = -\eps\sum_{i=1}^N\na\mu_{h,i}.
\end{align*}
Without regularization $\eps=0$, the total mass flux vanishes and $\bm{b}_{\rm skw}^\eps=0$. We define the discrete time derivative
\begin{align*}
  \pa_t^\tau w_h^{k+1} = \frac{1}{\tau}(w_h^{k+1}-w_h^k)
  \quad\mbox{for }w_h^{k+1},\,w_h^k\in\mathcal{V}_h
\end{align*}
and the nodal interpolation operator $I_h:C^0(\overline{\Omega})\to\mathcal{V}_h$.

Thus, our numerical scheme reads as
\begin{align}
  & \langle\pa_t^\tau \rho_{h,i}^{k+1},\psi_i\rangle
  - \langle \rho_{h,i}^k\bm{u}_h^{k+1},\na\psi_i\rangle
  + \sum_{j=1}^N\big\langle  M_{\eps,ij}(\pvec{\rho}_{h}^{*})\na\mu_{h,j}^{k+1},
  \na\psi_i\big\rangle = 0, \label{num.mass} \\
  & \langle \mu_{h,i}^{k+1},\xi_i\rangle 
  - \bigg\langle \frac{\pa f}{\pa\rho_i}(\pvec{\rho}_h^{k+1})
  + V_ip_h^{k+1},\xi_i\bigg\rangle = 0, \label{num.mu} \\
  & \frac{1}{2}\big\langle \partial_t^\tau\rho_h^{k+1}, I_h[\bm{u}_h^{k+1}
  \bm{v}]\big\rangle  + \big\langle \rho_h^k,I_h[
  \partial_t^\tau\bm{u}_h^{k+1}\bm{v}]\big\rangle \nonumber \\
  &\phantom{xx}+ \bm{b}_{\rm skw}(\rho_h^*\bm{u}_h^*;\bm{u}_h^{k+1},\bm{v})  
  + \bm{b}^\eps_{\rm skw}(\bm{J}_h^{k+1};\bm{u}_h^{k+1},\bm{v})
  - \langle  p_h^{k+1},\diver\bm{v}\rangle
  \label{num.mom} \\
  &\phantom{xx}
  + \big\langle \mathbb{S}(\vec{\rho}_{h}^{*},\na\bm{u}_h^{k+1}),
  \na\bm{v}\big\rangle + \sum_{i=1}^N\big\langle \rho_{h,i}^k\na(\mu_{h,i}^{k+1}-V_ip_h^{k+1}), 
  \bm{v}\big\rangle = 0, \nonumber \\
  &\langle \diver\bm{u}_h^{k+1},q\rangle
  = -\sum_{i,j=1}^N\big\langle  V_iM_{\eps,ij}(\pvec{\rho}_{h}^{*})
  \na\mu_{h,j}^{k+1},\na q\big\rangle \label{num.div}
\end{align}
for all $(\psi_i,\xi_i,\bm{v},q)\in\mathcal{V}_h\times\mathcal{V}_h
\times\mathcal{X}_h\times\mathcal{Q}_h$, $i=1,\ldots,N$. Here, we use the auxiliary variables
\begin{align*}
  \rho_h^* = \theta\rho_h^k + (1-\theta)\rho_h^{k+1}, \quad
  \bm{u}_h^* = \theta\bm{u}_h^k + (1-\theta)\bm{u}_h^{k+1} \quad
  \mbox{for }\theta\in[0,1].
\end{align*}
The case $\theta=1$ corresponds to an explicit discretization and $\theta=0$ to an implicit discretization. 


\subsection{Main results}

Let $\R_+:=[0,\infty)$. We impose the following assumptions.

\begin{enumerate}
\item[(A1)] Domain: $\Omega \subset \mathbb{R}^d$ ($d=2,3$) is a cube and identified with the $d$-dimensional torus. Functions on $\Omega$ are assumed to be periodic.

\item[(A2)] Diffusion matrix: The matrix $M:\mathbb{R}^N_+\to\mathbb{R}^{N\times N}$ with entries $M_{ij}\in C^1(\R_+^N)$ is symmetric and positive semi-definite with kernel spanned by $1=(1,\ldots,1)\in\mathbb{R}^N$, i.e.
\begin{equation*}
  \xi^\top M(\vec\rho)\xi \geq c(\vec{\rho})|\xi|^2\quad \mbox{with }c(\vec\rho)\ge 0 \quad\text{for all } 
  \xi\in\operatorname{span}\{1\}^\perp,
\end{equation*}
and $c(\vec{\rho})=0$ if and only if $\rho_i=0$ for at least one index $i=1,\ldots,N$.

\item[(A3)] Viscosity: The coefficients $\nu$, $\lambda:\mathbb{R}^N_+ \to \mathbb{R}_+$ are $C^1$-functions and positively bounded from above and below. The derivatives are uniformly bounded in $L^\infty(\R_+^N)$.
 
\item[(A4)] Specific volumes: The vector $\vec V=(V_1,\ldots,V_N)\in(0,\infty)^N$ is such that $\vec V\notin\operatorname{span}\{1\}$.

\item[(A5)] Internal energy density: $f\in C^{2,1}(\R_+^N)$ is homogeneous of degree one, i.e.\ $f(\vec\rho)=\sum_{i=1}^N \rho_i$ $(\pa f/\pa\rho_i)(\vec{\rho})$, and the Hessian $f''(\vec{\rho})$ is positive semidefinite with kernel spanned by the vector $\vec{\rho}$.

\item[(A6)] Initial data: $\rho_{i}^0\in C^0(\overline{\Omega})$ satisfies $0< c\leq \rho_{i}^0\leq V_i^{-1}$ for some $c>0$, $\sum_{i=1}^N V_i\rho_{i}^0=1$ a.e.\ in $\Omega$, and $\rho^0|\bm{u}^0|^2\in L^1(\Omega)$. 

\item[(A7)] Spatial discretization: $\mathcal{T}_h$  is quasi-uniform, i.e., there exists $\sigma>0$ such that $\sigma h\le\rho_K\le h_K\le h$ for $K\in\mathcal{T}_h$, where $\rho_K$ is the inner-circle radius, $h_K$ the diameter of $K\in\mathcal{T}_h$, and $h=\max_{K\in\mathcal{T}_h}h_K$.
\end{enumerate}

Recall that $I_h:C^0(\overline{\Omega})\to\mathcal{V}_h$ is the nodal interpolation operator and let $\pi_h^0:L^2(\Omega)\to\mathcal{V}_h$ be the $L^2$ projection in $\mathcal{V}_h$, defined by $\langle v-\pi_h^0v,w_h\rangle=0$ for all $v\in L^2(\Omega)$, $w_h\in\mathcal{V}_h$. The vector-valued $L^2$ projection on $\mathcal{X}_h$ is denoted by the same symbol. Our main result reads as follows.

\begin{theorem}[Error estimates]\label{thm.error}
Let (A1)--(A7) hold with $\eps> 0$. Let $(\vec{\rho}, \vec{\mu}, \bm{u}, p)$ be a regular solution to \eqref{1.ic}, \eqref{2.mass}--\eqref{2.div}, normalized such that $\langle p(t),1\rangle = 0$, such that
\begin{align*}
  & \vec{\rho},\, \bm{u}\in C^0([0,T];H^2(\Omega))\cap C^1([0,T];L^\infty(\Omega)) 
  \cap L^\infty(0,T;W^{1,\infty}(\Omega))\cap
  C^2([0,T];H^{-1}(\Omega)),  \\
  & \vec{\mu} \in C^0([0,T];H^2(\Omega)\cap W^{1,\infty}(\Omega)), \quad 
  p \in C^0([0,T];H^2(\Omega)). \notag
\end{align*}
Furthermore, let $h=c\tau$ with $\tau\leq \tau_0$ and $c>0$ and consider a family of discrete solutions $(\vec{\rho}_h,\vec{\mu}_h,\bm{u}_h,p_h)_{h>0}$ for the initial data $(\pvec{\rho}_h^0,\bm{u}^0_h) = (I_h\pvec{\rho}^0,\pi^0_h\bm{u}^0)$. Suppose that the lower bounds $\rho_{h,i}^k \geq c_i >0$ are uniform in $(h,\tau)$. Then there exists $C(\eps)>0$ such that for all $0<\tau\le\tau_0$ and $h=c\tau$,
\begin{align*}
  &\max_{k=1,\ldots,n}
  \big(\|\pvec{\rho}_h^k-\vec{\rho}(t^k)\|_{0,2}^2 
  + \|\bm{u}_h^k-\bm{u}(t^k)\|_{0,2}^2 \big) \\
  &\phantom{xx}+ \sum_{k=1}^{n}\tau\bigg(
  \|\vec{\mu}_h^k-\vec{\mu}(t^k)\|_{1,2}^2 
  + \|\bm{u}_h^k-\bm{u}(t^k)\|_{1,2}^2 
  + \|p_h^k-p(t^k)\|_{0,2}^2\bigg) \leq C(\eps)(h^2 + \tau^2).
  \end{align*}
\end{theorem}

\begin{remark}[Existence of discrete solutions]\label{rem.ex}\rm
Theorem \ref{thm.error} assumes the existence of discrete solutions with strictly positive densities. A local existence result can be shown by the implicit function theorem. Indeed, let $(\pvec{\rho}_h^k,\bm{u}_h^k)$ be given with $\rho_{h,i}^k>0$ in $\Omega$ for $i=1,\ldots,N$. The scheme \eqref{num.mass}--\eqref{num.div} can be formulated as $F(z,\tau)=0$, where $z$ is the nodal vector of $(\pvec{\rho}_h^{k+1}, \vec\mu_h^{k+1}, \bm{u}_h^{k+1}, p_h^{k+1})$ and $F$ is continuously differentiable in a neighborhood of positive densities. For $\tau=0$, a solution $z_0$ is given by $\pvec{\rho}_h^{k+1}=\pvec{\rho}_h^k$, $\bm{u}_h^{k+1} = \bm{u}_h^k$, where $p_h^{k+1}$ and $\vec\mu_h^{k+1}$ are determined by \eqref{num.mu} and \eqref{num.div}. It can be seen that the Jacobian of $F$ is invertible, and there exists $\tau_k>0$ such that the scheme possesses a solution for all $0<\tau<\tau_k$. By continuity, the densities $\rho_{h,i}^{k+1}$ remain positive. However, $\tau_k$ depends on $h$ and $\min_{i}\min_{x\in\Omega}\rho_{h,i}^k(x)$, and the discrete energy inequality does not provide a lower bound for the densities that is uniform in $k$. Consequently, the time step sizes may degenerate, and we cannot conclude the existence of a discrete solution on $[0,T]$. Topological degree arguments as in \cite{GGHL08} face the same difficulty, since they require a priori bounds that keep the densities strictly positive, uniformly with respect to the homotopy parameter. Positivity is obtained automatically in schemes formulated in the entropy variables \cite{BPS22}, but in the quasi-incompressible setting, such a formulation is not straightforward, since the chemical potentials contain the pressure, and the volume constraint \eqref{1.incom} has to be preserved. This explains the conditions imposed in Theorem \ref{thm.error}.
\end{remark}


\section{Properties of the numerical scheme}\label{sec.prop}

The numerical scheme corresponds to the scheme proposed in \cite{BJL25}, except for the $\eps$-regularization of the mobility matrix, the discretization of the convective term, and a lower polynomial order of the velocity space. The new discretization of the convective term with the total mass flux is important for the error analysis; see Section \ref{sec.drei}. Still, the same structure-preserving properties proved in \cite{BJL25} hold for the modified scheme.

\begin{theorem}[Discrete identities]\ 
Let Assumptions (A1)--(A6) hold and let $(\pvec{\rho}_{h}^{k}, \bm{u}_h^{k})\in$ $\Pi_\tau^1(\mathcal{V}_{h,+}^N \times\mathcal{X}_h)$, $(\vec\mu_{h}^k, p_h^k)\in  \Pi_\tau^0(\mathcal{V}_h^N \times \mathcal{Q}_h)$ be a discrete solution to \eqref{num.mass}--\eqref{num.div}. Then the following identities hold for $k\in \{0,\ldots,n-1\}$:
\begin{itemize}
\item Partial mass conservation: For all $i=1,\ldots,N$,
\begin{align*}
  \langle \rho_{h,i}^{k+1},\mathrm{1}\rangle
  = \langle \rho_{h,i}^{k},\mathrm{1}\rangle.
\end{align*}
\item Total mass conservation equation: For all $q\in\mathcal{V}_h$,
\begin{align*}
  \langle \partial_t^\tau\rho_h^{k+1},q\rangle
  - \langle \rho_h^k\bm{u}_h^{k+1}+\bm{J}_h^{k+1},\na q\rangle = 0.
\end{align*}
\item Constraint: For all $q\in\mathcal{Q}_h$,
\begin{align*}
 \bigg\langle \partial_t^\tau\bigg(\sum_{i=1}^N V_i\rho_{h,i}^{k+1} 
  - 1 \bigg),q\bigg\rangle 
  = \bigg\langle \bigg(\sum_{i=1}^N V_i\rho_{h,i}^{k} - 1\bigg)
  \bm{u}_h^{k+1},\na q\bigg\rangle.
\end{align*}
\item Discrete energy equality:
\begin{align*}
  \pa_t^\tau E_h^{k+1} &= -\big\langle \mathbb{S}(\pvec{\rho}_{h}^{*}, \na\bm{u}_h^{k+1}),\na\bm{u}_h^{k+1}\big\rangle 
  - \sum_{i,j=1}^N\big\langle  M_{\eps,ij}(\pvec{\rho}_{h}^{*})
  \na\mu_{h,i}^{k+1},\na\mu_{h,j}^{k+1}\big\rangle 
  - \mathcal{D}_{\rm num} \le 0,
\end{align*}
where the discrete energy and numerical dissipation are given by
\begin{align*}
  E_h^k &= \frac12\langle\rho_h^k,I_h[|\bm{u}_h^k|^2]\rangle
  + \langle f(\pvec{\rho}_{h}^{k}),1\rangle, \\
  \mathcal{D}_{\rm num} &= \frac{1}{2\tau}\big\langle\rho_h^k,
  I_h[|\bm{u}_h^{k+1}-\bm{u}_h^k|^2]\big\rangle
  + \frac{1}{2\tau}\sum_{i,j=1}^N
  \bigg\langle \frac{\pa^2 f}{\pa\rho_i\pa\rho_j}
  (\vec\omega_h)(\rho_{h,i}^{k+1}-\rho_{h,i}^k),
  \rho_{h,j}^{k+1}-\rho_{h,j}^k\bigg\rangle,
\end{align*}
and $\vec\omega_h$ is a convex combination of $\pvec{\rho}_{h}^{k}$ and $\pvec{\rho}_{h}^{k+1}$.
\end{itemize} 
\end{theorem}

\begin{theorem}[A priori bounds]\quad 
Let Assumptions (A1)--(A6) hold and  let $(\pvec{\rho}_{h}^{k}, \bm{u}_h^{k})\in$ $\Pi_\tau^1(\mathcal{V}_{h,+}^N\times\mathcal{X}_h)$, $(\vec\mu_{h}^k, p_h^k)\in \Pi_\tau^0(\mathcal{V}_h^N\times\mathcal{Q}_h)$ be a solution to \eqref{num.mass}--\eqref{num.div}. Then, for all $k\in \{0,\ldots,n\}$:
\begin{itemize}
\item Pointwise constraint:
\begin{align*}
  \sum_{i=1}^NV_i\rho_{h,i}^k(x)=1\quad\mbox{for all }x\in\Omega.
\end{align*}
\item Pointwise bounds:  
\begin{align*}
  \rho_{h,i}^k\le V_i^{-1}\quad\mbox{for }i=1,\ldots,N, \quad 
  V_{\rm max}^{-1}\le\rho_h^k\le V_{\rm min}^{-1},
\end{align*}
where $V_{\rm min}=\min\{V_1,\ldots,V_N\}$ and $V_{\rm max}=\max\{V_1,\ldots,V_N\}$. 
\item A priori estimates:
\begin{align*}
  \max_{k=1,\ldots,n}\big(\|\bm{u}_h^k\|_{0,2}^2
  + \|(\rho_{h}^k)^{1/2}\bm{u}_h^k\|_{0,2}^2\big)
  + \sum_{k=1}^{n}\tau\|\na\bm{u}_h^k\|_{0,2}^2
  \le C(E_h^0).
\end{align*}
\end{itemize}
\end{theorem}

\begin{proof}
	We follow the argumentation of \cite{BJL25}. In fact, we only need to prove the $L^2$ estimate for the velocity. For this, we use the equivalences
	\begin{equation*}
		\rho_{\min}\langle I_h[\bm{u}^2],1 \rangle \leq \langle \rho,I_h[\bm{u}^2] \rangle \leq \rho_{\max}\langle I_h[\bm{u}^2],1 \rangle, \quad 
		c\|\bm{u}\|_{0,2}^2 \leq \langle I_h[\bm{u}^2],1 \rangle \leq C\|\bm{u}\|_{0,2}^2
	\end{equation*}
    for some constants $\rho_{\min}$, $\rho_{\max}$, $c$, $C>0$. The first equivalence is a simple estimation, while the second one is a standard result for mass-lumped inner products, cf.\ \cite{GT24}. A combination of these inequalities leads to
	\begin{align}
		c\rho_{\min}\|\bm{u}\|_{0,2}^2 \leq \langle \rho,I_h[\bm{u}^2] \rangle \leq C\rho_{\max}\|\bm{u}\|_{0,2}^2, \label{eq:Interpequiv}
	\end{align}
which finishes the proof.
\end{proof}

\section{Error estimates}\label{sec.est1}

In this section, we prove Theorem \ref{thm.error}. For this, we split the total error into the projection error and the evolution error. First, we recall some estimates for the projection error. Then we define the perturbed system, satisfied by the projected solution, and compute the discrete relative energy equality involving the residuals. Next, we estimate the residual errors, which yields a discrete relative energy inequality; applying the discrete Gronwall lemma then gives the desired bounds. The final step of the proof is the computation of the pressure errors. 

\subsection{Projection error}

We estimate the difference of the solution $(\vec{\rho}, \vec{\mu}, \bm{u}, p)$ and $(\widehat\rho_{h,i}^k, \widehat{\bm{u}}_h^k,$ $\widehat{\mu}_{h,i}^k, \widehat{p}_h^k)$, where
\begin{align*}
  \widehat\rho_{h,i}^k = I_h\rho_i(t^k), \quad \widehat{\bm{u}}_h^k = \pi^0_h \bm{u}(t^k), \quad 
  \widehat\mu_{h,i}^k = \pi^0_h \mu_i(t^k), \quad 
  \widehat p_h^k = \pi^0_h p(t^k),
\end{align*}
recalling that $I_h:C^0(\overline{\Omega})\to\mathcal{V}_h$ denotes the nodal interpolation operator and $\pi_h^0:L^2(\Omega)\to\mathcal{V}_h$ is the $L^2$ projection. Restricting $\pi_h^0$ to zero-mean functions $u\in L^2(\Omega)$ with $\langle u,1\rangle=0$, the projection maps to $\mathcal{Q}_h$. The vector-valued $L^2$ projection for the velocity field is denoted by the same symbol $\pi_h^0$.

\begin{lemma}
Let $m^+=\max\{m,0\}$ for $m\in\{-1,0,1\}$ and $1\le p\le\infty$. Under the regularity assumption of Theorem \ref{thm.error} and Assumption (A7), the following projection errors hold for $i=1,\ldots,N$:
\begin{align}
  \|\widehat\rho_{h,i}^k-\rho_i(t^k)\|_{m,p}
  + \|\widehat{\bm{u}}_{h}^k-\bm{u}(t^k)\|_{m,p} &\le Ch^{2-m^+}, 
  \label{eq:projH1} \\
  \|\pa_t\widehat\rho_{h,i}^k
  - \pa_t\rho_i(t^k)\|_{-1,2} 
  + \|\pa_t\widehat{\bm{u}}_{h}^k
  - \pa_t\bm{u}(t^k)\|_{-1,2} &\le C(h^2\tau^{-1}+\tau), 
  \label{eq:projtime} \\
  \|\widehat\rho_{h,i}^*-\rho_i(t^k)\|_{0,p} 
  + \|\widehat{\bm{u}}_{h}^*-\bm{u}(t^k)\|_{0,p} &\le C(h^2+\tau), 
  \label{eq:projshift} \\
  \|\widehat\mu_{h,i}^k-\mu_i(t^k)\|_{m,p} 
  + \|\widehat p_{h}^k-p(t^k)\|_{m,p} &\le Ch^{2-m^+}.
  \label{eq:projH1mup}
\end{align}
Furthermore, the following stability estimates hold:
\begin{equation}\label{eq:projstab}
\begin{aligned}
  \|\widehat\rho_{h,i}^k\|_{1,\infty}
  &\leq C\|\rho_i(t^k)\|_{1,\infty}, \quad \|\widehat\mu_{h,i}^k\|_{1,\infty}
  \leq C\|\mu_i(t^k)\|_{1,\infty}, \\
  \|\widehat{\bm{u}}_{h}^k\|_{1,\infty}
  &\leq C\|\bm{u}(t^k)\|_{1,\infty}, \quad 
  \|\pa_t^\tau\widehat\rho_{h,i}^k\|_{0,\infty} 
  + \|\pa_t^\tau\widehat{\bm{u}}_{h}^k\|_{0,\infty} \leq C,
\end{aligned}
\end{equation}
where the constant $C>0$ does not depend on $h$ and $\tau$. 
\end{lemma}

The first three inequalities of the lemma follow from error estimates for the interpolation and projection operators (see, e.g., \cite{BrSc08}). The temporal estimates are based on a Taylor expansion and the interpolation estimates. The stability of $\pi_h^0$ in $W^{1,\infty}(\Omega)$ uses the quasi-uniformity assumption (A7).


\subsection{Perturbed system}

First, we recall the numerical scheme: Let $(\pvec{\rho}_{h}^0,\bm{u}_{h}^0) \in \mathcal{V}_h^{N} \times\mathcal{X}_h$ be given. Find $(\pvec{\rho}_h^{k+1},\bm{u}_h^{k+1}) \in \Pi^1_\tau(\mathcal{V}_h^{N}\times\mathcal{X}_h)$ and $(\vec{\mu}_{h}^{k+1},p_h^{k+1}) \in \Pi_\tau^0(\mathcal{V}_h^N \times\mathcal{Q}_h)$ such that 
\begin{align}
  & \langle \pa_t^\tau\rho_{h,i}^{k+1},\psi_i \rangle 
  - \langle \rho_{h,i}^{k}\bm{u}_h^{k+1}, \na\psi_i\rangle 
  + \textstyle\sum_{j=1}^N\langle  M_{\eps,ij}(\pvec{\rho}_{h}^*)
  \na\mu_{h,j}^{k+1},\na\psi_i\rangle = 0,  \label{N1} \\
  & \langle \mu_{h,i}^{k+1},\xi_i \rangle 
  - \langle (\pa f/\pa\rho_i)(\pvec{\rho}_{h}^{k+1})
  + V_ip_h^{k+1},\xi_i \rangle = 0, \label{N2} \\
  & \tfrac{1}{2}\langle \pa_t^\tau\rho_h^{k+1},I_h[\bm{u}_h^{k+1}
  \bm{v}]\rangle 
  + \langle \rho_h^k,I_h[\pa_t^\tau\bm{u}_h^{k+1}\bm{v}]\rangle
  + \bm{b}_{\rm skw}(\rho_h^*\bm{u}_h^*;\bm{u}_h^{k+1},\bm{v}) 
  + \bm{b}^\eps_{\rm skw}(\bm{J}_h^{k+1};\bm{u}_h^{k+1},\bm{v})
  \label{N3} \\
  &\phantom{xx}+ \langle \mathbb{S}(\pvec{\rho}_{h}^*,
  \na\bm{u}_h^{k+1}),\na\bm{v} \rangle 
  - \langle p_h^{k+1},\diver\bm{v} \rangle 
  + \textstyle\sum_{i=1}^N\langle 
  \rho_{h,i}^{k}\na(\mu_{h,i}^{k+1} - V_ip_h^{k+1}), 
  \bm{v}\rangle = 0, \nonumber \\
  & \langle \diver\bm{u}_h^{k+1},q \rangle 
  + \textstyle\sum_{i,j=1}^N\langle  
  V_iM_{\eps,ij}(\pvec{\rho}_{h}^*)\na\mu_{h,j}^{k+1},
  \na q \rangle = 0 \label{N4}
\end{align}
for all $((\psi_{i})_i,(\xi_i)_i,\bm{v},q)\in \mathcal{V}_h^N\times \mathcal{V}_h^N \times\mathcal{X}_h\times\mathcal{Q}_h$.

Next, let the perturbed initial data $(\pvec{\widehat{\rho}}_{h}^0, \widehat{\bm{u}}_{h}^0) \in \mathcal{V}_h^{N}\times\mathcal{X}_h$ and the projected solutions $(\pvec{\widehat\rho}_{h}^{k+1},\widehat{\bm{u}}_h^{k+1})\in \Pi^1_\tau(\mathcal{V}_h^{N} \times\mathcal{X}_h)$ and $(\pvec{\widehat\mu}_h^{k+1},\widehat p_h^{k+1})\in \Pi^0_\tau(\mathcal{V}_h^N\times\mathcal{Q}_h)$ be given. Inserting this solution into the equations yields some errors, which we denote by $r_{1,i}^{k+1}$, $r_{2,i}^{k+1}$, $\bm{r}_3^{k+1}$, and $r_4^{k+1}$ (see \eqref{5.res} below for the precise definitions):
\begin{align}
  & \langle\pa_t^\tau\widehat\rho_{h,i}^{k+1},\psi_i \rangle 
  - \langle \rho_{h,i}^{k}\widehat{\bm{u}}_h^{k+1}, 
  \na\psi_i\rangle + \textstyle\sum_{j=1}^N\big\langle  M_{\eps,ij}(\pvec{\widehat{\rho}}_{h}^*)\na\widehat\mu_{h,j}^{k+1},
  \na\psi_i\big\rangle \label{pertN1} \\ 
  &\phantom{xx}= \langle r_{1,i}^{k+1}+ s_{1,i}^{k+1}, 
  \psi_i\rangle, \nonumber \\
  &\langle\widehat\mu_{h,i}^{k+1},\xi_i \rangle 
  - \big\langle(\pa f/\pa\rho_i)(\pvec{\widehat\rho}_{h}^{k+1}) 
  \red{+} V_i\widehat{p}_h^{k+1},\xi_i \big\rangle 
  = \langle r_{2,i}^{k+1}, \xi_i\rangle, \label{pertN2} \\
  &  \tfrac{1}{2}\langle\pa_t^\tau\rho_h^{k+1},
  I_h[\widehat{\bm{u}}_h^{k+1}\bm{v}]\rangle 
  + \langle\rho_h^k,
  I_h[\pa_t^\tau\widehat{\bm{u}}_h^{k+1}\bm{v}]\rangle
  + \bm{b}_{\rm skw}(\rho_h^*\bm{u}_h^*;
  \widehat{\bm{u}}_h^{k+1},\bm{v}) 
  \label{pertN3} \\
  &\phantom{xx} + \bm{b}^\eps_{\rm skw}(\bm{J}_h^{k+1}; \widehat{\bm{u}}_h^{k+1},\bm{v})
  + \langle \mathbb{S}(\pvec{\rho}_{h}^*,
  \na\widehat{\bm{u}}_h^{k+1}),\na\bm{v} \rangle 
  - \langle\widehat{p}_h^{k+1},\diver\bm{v}\rangle \nonumber \\
  &\phantom{xx}+ \textstyle\sum_{i=1}^N\langle 
  \rho_{h,i}^{k}\na(\widehat\mu_{h,i}^{k+1} 
  - V_i\widehat{p}_h^{k+1}), \bm{v}\rangle 
  = \langle \bm{r}_{3}^{k+1} + \bm{s}_{3}^{k+1}, 
  \bm{v}\rangle, \nonumber \\
  & \langle \diver\widehat{\bm{u}}_h^{k+1},q \rangle 
  + \textstyle\sum_{i,j=1}^N
  \big\langle V_iM_{\eps,ij}(\pvec{\widehat\rho}_{h}^*)
  \na\widehat\mu_{h,j}^{k+1},\nabla q \big\rangle 
  = \langle r_{4}^{k+1}, q\rangle \label{pertN4}
\end{align}
for all $((\psi_{i})_i,(\xi_i)_i,\bm{v},q)\in \mathcal{V}_h^N\times \mathcal{V}_h^N \times\mathcal{X}_h\times\mathcal{Q}_h$, where the variables $s_{1,i}^{k+1}$ and $\bm{s}_3^{k+1}$ are defined by
\begin{align}
  \langle s_{1,i}^{k+1},\psi_i \rangle 
  &= -\langle (\rho_{h,i}^{k}-\widehat\rho_{h,i}^{k})
  \widehat{\bm{u}}_h^{k+1}, \na\psi_i \rangle, \label{s1} \\
  \langle \bm{s}_3^{k+1},\bm{v}\rangle 
  &= \tfrac{1}{2}\langle \pa_t^\tau(\rho_h^{k+1}-\widehat\rho_h^{k+1}),
  I_h[\widehat{\bm{u}}_h^{k+1}\bm{v}]\rangle 
  + \langle \rho_h^k-\widehat\rho_h^k,I_h[\pa_t^\tau\widehat{\bm{u}}_h^{k+1}\bm{v}]\rangle\label{s3} \\
  &\phantom{xx}+ \bm{b}_{\rm skw}(\rho_h^*\bm{u}_h^* 
  - \widehat\rho_h^* \widehat{\bm{u}}_h^*;
  \widehat{\bm{u}}_h^{k+1},\bm{v})
  + \bm{b}^\eps_{\rm skw}(\bm{J}_h^{k+1}
  - \widehat{\bm{J}}_h^{k+1};\widehat{\bm{u}}_h^{k+1},\bm{v}) 
  \nonumber \\
  &\phantom{xx}+ \big\langle\mathbb{S}(\vec{\rho}_{h}^*,
  \na\widehat{\bm{u}}_h^{k+1}) 
  - \mathbb{S}(\vec{\widehat{\rho}}_{h}^*,
  \na\widehat{\bm{u}}_h^{k+1}),\na\bm{v}\big\rangle \nonumber \\
  &\phantom{xx}+ \textstyle\sum_{i=1}^N\big\langle (\rho_{h,i}^{k}-\widehat\rho_{h,i}^k)
  \na(\widehat\mu_{h,i}^{k+1} - V_i\widehat p_h^{k+1}),\bm{v}\big\rangle.
  \nonumber 
\end{align}

Since $\widehat\rho_h^k$ is the nodal interpolation of the exact solution, the volume-filling constraint is preserved at the nodal values. As any point is a convex combination of the nodal points, the constraint holds pointwise everywhere. Choosing the test function $\psi_i=V_iq$ in \eqref{pertN1} and taking into account \eqref{pertN4} leads to
\begin{align}\label{r1r4}
  \langle r_4^{k+1},q\rangle 
  = \sum_{i=1}^NV_i\langle(r_{1,i}^{k+1}+s_{1,i}^{k+1}),q\rangle
  \quad\mbox{for all }q\in\mathcal{Q}_h.
\end{align}
This identity will be used later.


\subsection{Relative energy}

We introduce the relative energy
\begin{align*}
  \mathcal{E}\big(\vec{\rho},\bm{u}|\vec{\widehat{\rho}},
  \widehat{\bm{u}}\big) := \bigg\langle f(\vec{\rho}) 
  - f(\vec{\widehat{\rho}}) 
  - \sum_{i=1}^N\frac{\pa f(\vec{\widehat{\rho}})}{\pa\rho_i}
  (\rho_i-\widehat{\rho}_i),1\bigg\rangle 
  + \frac12\big\langle\rho,I_h[|\bm{u}-\widehat{\bm{u}}|^2]\big\rangle.
\end{align*}
We claim that for given $\vec\rho$ and $\vec{\widehat{\rho}}$ such that $0<c\le \rho_i,\,\widehat\rho_i\le V_i^{-1}$, there exist constants $c_1>0$, $C_1>0$ such that
\begin{align}\label{relenerL2}
  c_1\big(\|\vec\rho-\vec{\widehat{\rho}}\|_{0,2}^2
  + \|\bm{u}-\widehat{\bm{u}}\|_{0,2}^2\big)
  \le \mathcal{E}\big(\vec{\rho},\bm{u}|\vec{\widehat{\rho}},
  \widehat{\bm{u}}\big)
  \le C_1\big(\|\vec\rho-\vec{\widehat{\rho}}\|_{0,2}^2
  + \|\bm{u}-\widehat{\bm{u}}\|_{0,2}^2\big).
\end{align}

For the relative kinetic energy, this is obvious using \eqref{eq:Interpequiv}. For the internal mixing energy we conclude from Taylor expansion that
\begin{align}\label{eq:relhessian}
  \bigg\langle f(\vec\rho)-f(\vec{\widehat{\rho}})
  - \sum_{i=1}^N\frac{\pa f}{\pa\rho_i}(\vec{\widehat\rho})(\rho_i-\widehat\rho_i),
  1\bigg\rangle = \frac12\big\langle f''(\vec\omega)(\vec\rho
  - \vec{\widehat\rho}),\vec\rho - \vec{\widehat\rho}\big\rangle
\end{align}
for some $\vec\omega$, where $f''(\vec\omega)$ denotes the Hessian matrix. The homogeneity of degree one of $f$ implies that
\begin{align*}
  0 = \frac{\pa}{\pa\rho_i}\bigg(\sum_{j=1}^N\rho_j
  \frac{\pa f}{\pa\rho_j}(\vec\rho)\bigg) 
  - \frac{\pa f}{\pa\rho_i}(\vec\rho)
  = \sum_{j=1}^N\rho_j\frac{\pa^2 f}{\pa\rho_j\pa\rho_i}(\vec\rho).
\end{align*}
Consequently, $\vec\rho$ is an element of the kernel of $f''(\vec\rho)$ and by assumption, the kernel is one-dimensional. Thus, \eqref{eq:relhessian} is bounded from below by the $L^2(\Omega)$ norm if $\vec\omega\neq c(\vec\rho-\vec{\widehat\rho})$ for any $c\in\R$. By contradiction, we deduce from  $\vec\omega=\theta\vec\rho+(1-\theta)\vec{\widehat\rho}$ for some $\theta\in[0,1]$ that
\begin{align*}
  0 = \vec\omega - c(\vec\rho-\vec{\widehat\rho})
  = (\theta-c)\vec\rho + (1-\theta+c)\vec{\widehat\rho}.
\end{align*}
After multiplication with $\vec{V}$ and using the constraint, this yields
\begin{align*}
  0 = (\theta-c)\vec\rho\cdot\vec V 
  + (1-\theta+c)\vec{\widehat\rho}\cdot\vec V
  = (\theta-c) + (1-\theta+c) = 1,
\end{align*}
which is a contradiction. We conclude that the Hessian can be estimated from above and below, using the uniform bounds. This shows the claim.


\subsection{Discrete relative energy identity}

The difference of the discrete relative energies at times $t^{k+1}$ and $t^k$ is given by
\begin{align*}
 \frac{1}{\tau}\mathcal{E}(\vec{\rho}_h,\bm{u}_h|
  \vec{\widehat\rho}_h,\widehat{\bm{u}}_h)\Big|_{t^k}^{t^{k+1}} 
  &= \frac{1}{2\tau}\bigg\langle\rho_h^{k+1}I_h[|\bm{u}_h^{k+1}
  -\widehat{\bm{u}}_h^{k+1}|^2] 
  - \rho_h^{k}I_h[|\bm{u}_h^{k}-\widehat{\bm{u}}_h^{k}|^2],1
  \bigg\rangle \\
  &\phantom{xx}+ \frac{1}{\tau}\bigg\langle f(\pvec{\rho}_h^{k+1}) 
  - f(\pvec{\widehat\rho}_h^{k+1}) 
  - \sum_{i=1}^N\frac{\pa f}{\pa\rho_i}(\pvec{\widehat\rho}_h^{k+1})
  (\rho_{h,i}^{k+1}-\widehat\rho_{h,i}^{k+1}), 1\bigg\rangle \nonumber \\
  &\phantom{xx}- \frac{1}{\tau}\bigg\langle 
  f(\pvec{\rho}_h^{k})-f(\pvec{\widehat\rho}_h^{k}) 
  - \sum_{i=1}^N\frac{\pa f}{\pa\rho_i}(\pvec{\widehat\rho}_h^{k})
  (\rho_{h,i}^k-\widehat\rho_{h,i}^k), 1\bigg\rangle \nonumber \\
  &=: I_1 + I_2 + I_3. \nonumber 
\end{align*}
The first term $I_1$ can be rewritten as $I_1=I_{11}+I_{12}$, where
\begin{align*}
  I_{11} &= \frac{1}{2}\big\langle\pa_t^\tau\rho_h^{k+1},I_h[ (\bm{u}_h^{k+1}-\widehat{\bm{u}}_h^{k+1})
  (\bm{u}_h^{k+1}-\widehat{\bm{u}}_h^{k+1})]\big\rangle \\
  &\phantom{xx}+ \big\langle \rho^{k}_h, I_h[\pa_t^\tau(\bm{u}_h^{k+1} 
  - \widehat{\bm{u}}_h^{k+1})\cdot
  (\bm{u}_h^{k+1}-\widehat{\bm{u}}_h^{k+1})] \big\rangle, \\
  I_{12} &= -\frac{\tau}{2}\big\langle \rho_h^k,I_h[|\pa_t^\tau(\bm{u}_h^{k+1}
  -\widehat{\bm{u}}_h^{k+1})|^2]\big\rangle.
\end{align*}
The last term corresponds to the (negative) numerical dissipation, so that we only need to reformulate the term $I_{11}$. We insert the test function $\bm{v}=\bm{u}_h^{k+1} - \widehat{\bm{u}}_h^{k+1}$ in \eqref{N3} and \eqref{pertN3} to replace the discrete time derivatives:
\begin{align*}
  I_{11} &= -\bm{b}_{\rm skw}(\rho_h^*\bm{u}_h^*;\bm{u}_h^{k+1} 
  - \widehat{\bm{u}}_h^{k+1},
  \bm{u}_h^{k+1}-\widehat{\bm{u}}_h^{k+1}) 
  - \bm{b}^\eps_{\rm skw}(\bm{J}_h^{k+1};\bm{u}_h^{k+1}-\widehat{\bm{u}}_h^{k+1},
  \bm{u}_h^{k+1}-\widehat{\bm{u}}_h^{k+1}) \\
  &\phantom{xx}- \langle \mathbb{S}(\pvec{\rho}_h^*,\na(\bm{u}_h^{k+1}
  - \widehat{\bm{u}}_h^{k+1})),
  \na(\bm{u}_h^{k+1}-\widehat{\bm{u}}_h^{k+1})\rangle 
  + \langle p_h^{k+1}-\widehat p_h^{k+1},\diver(\bm{u}_h^{k+1}-\widehat{\bm{u}}_h^{k+1}) \rangle \\
  &\phantom{xx}- \sum_{i=1}^N \big\langle \rho_{h,i}^k
  \na\big(\mu_{h,i}^{k+1} - \widehat\mu_{h,i}^{k+1} 
  - V_i(p_h^{k+1}-\widehat p_h^{k+1})\big),
  \bm{u}_h^{k+1}-\widehat{\bm{u}}_h^{k+1} \big\rangle \\ 
  &\phantom{xx}- \langle \bm{r}_3^{k+1}+\bm{s}_3^{k+1},
  \bm{u}_h^{k+1}-\widehat{\bm{u}}_h^{k+1} \rangle.
\end{align*}
Because of the skew-symmetry of $\bm{b}_{\rm skw}$ and $\bm{b}^\eps_{\rm skw}$, the first two terms on the right-hand side vanish. The fourth term can be reformulated by using the test function $q=p_h^{k+1}-\widehat{p}_h^{k+1}$ in \eqref{N4} and \eqref{pertN4}, leading to
\begin{align*}
  I_{11} &= - \langle \mathbb{S}(\pvec{\rho}_h^*,\na(\bm{u}_h^{k+1}
  - \widehat{\bm{u}}_h^{k+1})),
  \na(\bm{u}_h^{k+1}-\widehat{\bm{u}}_h^{k+1})\rangle \\
  &\phantom{xx}- \sum_{i,j=1}^N \big\langle M_{\eps,ij}(\pvec{\rho}_h^*)
  \na\mu_{h,j}^{k+1} - M_{\eps,ij}(\pvec{\widehat\rho}_h^*)
  \na\widehat\mu_{h,j}^{k+1}, 
  V_i\na(p_h^{k+1}-\widehat p_h^{k+1})\big\rangle \\
  &\phantom{xx}- \sum_{i=1}^N \langle \rho_{h,i}^k
  \na(\mu_{h,i}^{k+1}-\widehat\mu_{h,i}^{k+1}),
  \bm{u}_h^{k+1}-\widehat{\bm{u}}_h^{k+1} \rangle \\
  &\phantom{xx}+ \sum_{i=1}^N \langle \rho_{h,i}^k 
  V_i\na(p_h^{k+1}-\widehat p_h^{k+1}),
  \bm{u}_h^{k+1}-\widehat{\bm{u}}_h^{k+1} \rangle \\
  &\phantom{xx}- \langle \bm{r}_3^{k+1}+\bm{s}_3^{k+1},
  \bm{u}_h^{k+1}-\widehat{\bm{u}}_h^{k+1} \rangle 
  - \langle r_4^{k+1},p_h^{k+1}-\widehat p_h^{k+1} \rangle.
\end{align*}
We insert the test function $\psi_i=V_i(p_h^{k+1}-\widehat{p}_h^{k+1})$ in \eqref{N1} and \eqref{pertN1}; then the fourth term on the right-hand side cancels out, and we obtain
\begin{align}\label{4.I1}
  I_{11} &= - \langle \mathbb{S}(\pvec{\rho}_h^*,\na(\bm{u}_h^{k+1}
  - \widehat{\bm{u}}_h^{k+1})),
  \na(\bm{u}_h^{k+1}-\widehat{\bm{u}}_h^{k+1})\rangle \\
  &\phantom{xx}+ \sum_{i=1}^N \langle \pa_t^\tau(\rho_{h,i}^{k+1}-\widehat\rho_{h,i}^{k+1}),
  V_i(p_h^{k+1}-\widehat p_h^{k+1}) \rangle \nonumber \\
  &\phantom{xx}- \sum_{i=1}^N \langle \rho_{h,i}^k
  \na (\mu_{h,i}^{k+1}-\widehat\mu_{h,i}^{k+1}), \bm{u}_h^{k+1}-\widehat{\bm{u}}_h^{k+1} \rangle 
  + \mathcal{R}_1, \nonumber 
\end{align}
where the remainder $\mathcal{R}_1$ is defined by
\begin{align*}
  \mathcal{R}_1 &:= -\langle \bm{r}_3^{k+1}+\bm{s}_3^{k+1},
  \bm{u}_h^{k+1}-\widehat{\bm{u}}_h^{k+1} \rangle 
  - \langle r_4^{k+1}, p_h^{k+1}-\widehat p_h^{k+1} \rangle \\
  &\phantom{xx}+ \sum_{i=1}^N\langle r_{1,i}^{k+1}+s_{1,i}^{k+1},V_i(p_h^{k+1}
  - \widehat p_h^{k+1}) \rangle.
\end{align*}
The term $I_2+I_3$ is formulated as
\begin{align*}
  I_2+I_3 = \sum_{i=1}^N\bigg\langle \pa_t^\tau(\rho_{h,i}^{k+1}-\widehat\rho_{h,i}^{k+1}),
  \frac{\pa f}{\pa\rho_i}(\pvec{\rho}_h^{k+1}) 
  - \frac{\pa f}{\pa\rho_i}(\pvec{\widehat\rho}_h^{k+1})\bigg\rangle 
  + \mathcal{R}_2,
\end{align*}
where the remainder $\mathcal{R}_2$ is defined by
\begin{align*}
  \mathcal{R}_2 &= \frac{1}{\tau}\big\langle f(\pvec{\rho}_h^{k+1})-f(\pvec{\widehat\rho}_h^{k+1})
  -\sum_{i=1}^N\frac{\pa f}{\pa\rho_i}(\pvec{\widehat\rho}_h^{k+1})(\rho_{h,i}^{k+1}-\widehat\rho_{h,i}^{k+1}),1 \big\rangle \\
  &\phantom{xx}- \frac{1}{\tau}\big\langle
  f(\pvec{\rho}_h^{k})-f(\pvec{\widehat\rho}_h^{k})
  -\sum_{i=1}^N\frac{\pa f}{\pa\rho_i}(\pvec{\widehat\rho}_h^{k})(\rho_{h,i}^{k}-\widehat\rho_{h,i}^{k}),1 \big\rangle \\
  &\phantom{xx}- \frac{1}{\tau}\sum_{i=1}^N\bigg\langle (\rho_{h,i}^{k+1}-\widehat\rho_{h,i}^{k+1})
  - (\rho_{h,i}^{k}-\widehat\rho_{h,i}^{k}),
  \frac{\pa f}{\pa\rho_i}(\pvec{\rho}_h^{k+1}) 
  - \frac{\pa f}{\pa\rho_i}(\pvec{\widehat\rho}_h^{k+1}) \bigg\rangle. 
\end{align*}
Together with $I_{11}$ from \eqref{4.I1}, this gives
\begin{align*}
  &\frac{1}{\tau}\mathcal{E}(\vec{\rho}_h,\bm{u}_h|
  \vec{\widehat\rho}_h,\widehat{\bm{u}}_h)\Big|_{t^k}^{t^{k+1}}
  = I_1+I_2+I_3 = -\langle \mathbb{S}(\pvec{\rho}_h^{*}, \na(\bm{u}_h^{k+1}-\widehat{\bm{u}}_h^{k+1})), \na(\bm{u}_h^{k+1}-\widehat{\bm{u}}_h^{k+1})\rangle \\
  &\phantom{xx}+ \sum_{i=1}^N \bigg\langle \pa_t^\tau(\rho_{h,i}^{k+1}-\widehat\rho_{h,i}^{k+1}), 
  \frac{\pa f}{\pa\rho_i}(\pvec{\rho}_h^{k+1}) 
  - \frac{\pa f}{\pa\rho_i}(\pvec{\widehat\rho}_h^{k+1}) 
  + V_i(p_h^{k+1}-\widehat p_h^{k+1}) \bigg\rangle \\
  &\phantom{xx}- \sum_{i=1}^N \langle \rho_{h,i}^k
  \na (\mu_{h,i}^{k+1}-\widehat\mu_{h,i}^{k+1}), \bm{u}_h^{k+1}-\widehat{\bm{u}}_h^{k+1} \rangle 
  + I_{12} + \mathcal{R}_1 + \mathcal{R}_2.
\end{align*}
The second term on the right-hand side can be replaced after inserting the test function $\xi_i=\pa_t^\tau(\rho_{h,i}^{k+1}-\widehat\rho_{h,i}^{k+1})$ in \eqref{N2} and \eqref{pertN2}, and the third term is replaced after choosing $\psi_i=\mu_{h,i}^{k+1}-\widehat{\mu}_{h,i}^{k+1}+r_{2,i}^{k+1}$ in \eqref{N1} and \eqref{pertN1}. Taking into account relation \eqref{r1r4} between $r_{4}^{k+1}$ and $r_{1,i}^{k+1}+s_{1,i}^{k+1}$, it follows after a computation and multiplication with $\tau$ that
\begin{align}
  \mathcal{E}(\vec{\rho}_h,\bm{u}_h|
  \vec{\widehat\rho}_h,\widehat{\bm{u}}_h)\Big|_{t^k}^{t^{k+1}}
  &+ \tau\mathcal{D}(\vec{\mu}_h^{k+1}-\pvec{\widehat{\mu}}_h^{k+1},
  \bm{u}_h^{k+1}-\widehat{\bm{u}}_h^{k+1}) \label{discre11} \\
  &+ \frac{\tau^2}{2}\big\langle\rho_h^k,I_h[|\pa_t^\tau
  (\bm{u}_h^{k+1}-\widehat{\bm{u}}_h^{k+1})|^2]\big\rangle
  = J_1+\cdots+J_6, \nonumber 
\end{align}
where the energy dissipation is given by
\begin{align*}
  \mathcal{D}\big(\vec{\mu}_h^{k+1}-\pvec{\widehat{\mu}}_h^{k+1},
  \bm{u}_h^{k+1}-\widehat{\bm{u}}_h^{k+1}\big)
  &= \langle \mathbb{S}(\pvec{\rho}_h^{*}, \na(\bm{u}_h^{k+1}-\widehat{\bm{u}}_h^{k+1})),
  \na(\bm{u}_h^{k+1}-\widehat{\bm{u}}_h^{k+1})\rangle \\
  &\phantom{xx}+ \sum_{i,j=1}^N \langle 
  M_{\eps,ij}(\pvec{\rho}_h^*)
  \na(\mu_{h,j}^{k+1}-\widehat\mu_{h,j}^{k+1}),
  \na(\mu_{h,i}^{k+1}-\widehat\mu_{h,i}^{k+1}) \rangle,
\end{align*}
and the terms on the right-hand side of \eqref{discre11} read as
\begin{align*}
  J_1 &= -\tau\sum_{i,j=1}^N \langle M_{\eps,ij}(\pvec{\rho}_h^*)
  \na(\mu_{h,j}^{k+1}-\widehat\mu_{h,j}^{k+1}),
  \na r_{2,i}^{k+1} \rangle, \\
  J_2 &= \red{-}\tau\sum_{i,j=1}^N \langle 
  (M_{\eps,ij}(\pvec{\rho}_h^*)- M_{\eps,ij}(\pvec{\widehat\rho}_h^*))
  \na\widehat\mu_{h,j}^{k+1},\na(\mu_{h,i}^{k+1} 
  - \widehat\mu_{h,i}^{k+1} + r_{2,i}^{k+1}) \rangle, \\
  J_3 &= \tau\sum_{i=1}^N \langle \rho_{h,i}^k(\bm{u}_h^{k+1}
  - \widehat{\bm{u}}_h^{k+1}),\na r_{2,i}^{k+1} \rangle, \\
  J_4 &= \tau\mathcal{R}_2, \\
  J_5 &= -\tau\langle \bm{r}_3^{k+1} + \bm{s}_3^{k+1},\bm{u}_h^{k+1}
  - \widehat{\bm{u}}_h^{k+1} \rangle, \\
  J_6 &= -\tau\sum_{i=1}^N\langle r_{1,i}^{k+1}+s_{1,i}^{k+1},\mu_{h,i}^{k+1}
  -\widehat\mu_{h,i}^{k+1}+r_{2,i}^{k+1} \rangle.
\end{align*}
For the following estimates, we will often use the following bounds, which follow from the uniform lower and upper bounds for the densities:
\begin{align}\label{4.D}
  C_1\big(\|\na&(\bm{u}_h^{k+1}-\widehat{\bm{u}}_h^{k+1})\|_{0,2}^2
  + \|\na(\vec{\mu}_h^{k+1}-\pvec{\widehat{\mu}}_h^{k+1})\|_{0,2}^2
  \big) \\
  &\le \mathcal{D}\big(\vec{\mu}_h^{k+1}-\pvec{\widehat{\mu}}_h^{k+1},
  \bm{u}_h^{k+1}-\widehat{\bm{u}}_h^{k+1}\big) \nonumber \\
  &\le C_2\big(\|\na(\bm{u}_h^{k+1}-\widehat{\bm{u}}_h^{k+1})\|_{0,2}^2
  + \|\na(\vec{\mu}_h^{k+1}-\pvec{\widehat{\mu}}_h^{k+1})\|_{0,2}^2\big). 
  \nonumber 
\end{align}


\subsection{Discrete relative energy inequality}\label{sec.drei}

We estimate the terms $J_1,\ldots,J_6$ in \eqref{discre11}. We infer from Young's inequality for some $\delta>0$ that
\begin{align*}
  J_1 \le \delta\tau \mathcal{D}
  \big(\vec{\mu}_h^{k+1}-\pvec{\widehat\mu}_h^{k+1},0\big)
  + C(\delta)\tau\sum_{i,j=1}^N\langle M_{\eps,ij}(\pvec{\rho}_h^*)
  \na r_{2,i}^{k+1},\na r_{2,j}^{k+1} \rangle.
\end{align*}
For the term $J_2$, we observe that, by construction,
$M_{\eps,ij}(\pvec{\rho}_h^*) - M_{\eps,ij}(\pvec{\widehat\rho}_h^*) 
= M_{ij}(\pvec{\rho}_h^*)-M_{ij}(\pvec{\widehat\rho}_h^*)$, which yields
\begin{align*}
  J_2 \le \tau\big\|\big(M(\pvec{\rho}_h^*) 
  - M(\pvec{\widehat\rho}_h^*)\big)
  \na\pvec{\widehat\mu}_h^{k+1}\cdot
  \na\big(\vec{\mu}_h^{k+1}-\pvec{\widehat\mu}_h^{k+1}
  + r_{2}^{k+1}\big)\big\|_{0,1}.
\end{align*}
We estimate $J_2$ similarly as in \cite[page 19]{Dru23}. Let $M^\dag$ be the pseudo-inverse of $M$. Then, by definition, $MM^\dag M=M$, and consequently, $\langle(M-\widehat{M})\na a,\na b\rangle = \langle M^\dag(M-\widehat{M})\na a,M\na b\rangle$. We infer from the Cauchy--Schwarz and Young inequalities with $\delta>0$ as well as estimate \eqref{4.D} that
\begin{align*}
  J_2 &\le \tau\big\|M^\dag(\pvec{\rho}_h^*)
  \big(M(\pvec{\rho}_h^*) - M(\pvec{\widehat\rho}_h^*)\big)
  \na\pvec{\widehat\mu}_h^{k+1}\cdot M(\vec{\rho}_h^*)
  \na\big(\vec{\mu}_h^{k+1}-\pvec{\widehat\mu}_h^{k+1}
  + r_{2}^{k+1}\big)\big\|_{0,1} \\
  &\le \delta\tau\mathcal{D}
  \big(\vec{\mu}_h^{k+1}-\pvec{\widehat\mu}_h^{k+1},0\big) 
  + C(\delta)\tau\|M^\dag(\pvec{\rho}_h^*)\|_{0,\infty}^2\big\|M(\pvec{\rho}_h^*) - M(\pvec{\widehat\rho}_h^*)\big\|_{0,2}^2
  \|\na \pvec{\widehat\mu}_h^{k+1}\|_{0,\infty}^2 \\
  &\phantom{xx}+ C(\delta)\tau\|\na r_2^{k+1}\|_{0,2}^2 \\
  &\le \delta\tau\mathcal{D}
  \big(\vec{\mu}_h^{k+1} - \pvec{\widehat\mu}_h^{k+1},0\big) 
  + C(\delta)\tau\|
  \na\pvec{\widehat\mu}_h^{k+1}\|_{0,\infty}^2
  \big(\mathcal{E}\big(\pvec{\rho}_h^{k+1},\bm{u}_h^{k+1}|
  \pvec{\widehat\rho}_h^{k+1},\widehat{\bm{u}}_h^{k+1}\big)
  + \mathcal{E}\big(\pvec{\rho}_h^{k},\bm{u}_h^{k}|
  \pvec{\widehat\rho}_h^{k},\widehat{\bm{u}}_h^{k}\big)\big) \\
  &\phantom{xx}+ C(\delta)\tau\|\na r_2^{k+1}\|_{0,2}^2,
\end{align*}
where the last step follows from the uniform bounds for the densities, the Lipschitz continuity of $M$, and relations \eqref{relenerL2} between the relative energy and the $L^2(\Omega)$ difference of the densities.  The term $J_3$ can be estimated as
\begin{align*}
  J_3 \le \sum_{i=1}^N\tau\big(\|r_{2,i}^{k+1}\|_{1,2}^2 
  + C\|\rho_{h,i}^k\|_{0,\infty}^2
  \|\bm{u}_h^{k+1} - \widehat{\bm{u}}_h^{k+1}\|_{0,2}^2\big).
\end{align*}
The term $J_4$ can be seen as the discrete time derivative of the relative internal energy. Then, similarly as in \cite[Lemma 4.4]{EGKP23}, 
\begin{align*}
  J_4 &\le \sum_{i=1}^N\bigg(\bigg\langle 
  \frac{\pa f}{\pa\rho_i}(\pvec{\rho}_h^{k+1}) 
  - \frac{\pa f}{\pa\rho_i}(\pvec{\widehat\rho}_h^{k+1}) 
  - \sum_{j=1}^N \frac{\pa^2 f}{\pa\rho_i\pa\rho_j}
  (\pvec{\widehat\rho}_h^{k+1})
  (\rho_{h,j}^{k+1} - \widehat\rho_{h,j}^{k+1}), 
  \widehat\rho_{h,i}^{k+1} - \widehat\rho_{h,i}^k\bigg\rangle \\
 &\phantom{xx} + C\|\widehat\rho_{h,i}^{k+1} 
  - \widehat\rho_{h,i}^k\|_{0,\infty}
  \Big(\mathcal{E}(\pvec{\rho}_h^{k+1},0|\pvec{\widehat\rho}_h^{k},0) 
  + \mathcal{E}(\pvec{\rho}_h^{k},0|\pvec{\widehat\rho}_h^{k},0)  
  + C\|\widehat\rho_{h,i}^{k+1}-\widehat\rho_{h,i}^k\|_{0,2}^2\Big)
  \bigg) \\
  & \leq C\tau\sum_{i=1}^N
  \|\pa_t^\tau\widehat\rho_{h,i}^{k+1}\|_{0,\infty}
  \big(\mathcal{E}(\pvec{\rho}_h^{k+1},0|\pvec{\widehat\rho}_h^{k},0) 
  + \mathcal{E}(\pvec{\rho}_h^{k},0|\pvec{\widehat\rho}_h^{k},0)\big) \\
  &\phantom{xx}+ C\tau^3\sum_{i=1}^N
  \|\pa_t^\tau\widehat\rho_{h,i}^{k+1}\|_{0,\infty}
  \|\pa_t^\tau\widehat\rho_{h,i}^{k+1}\|_{0,2}^2.
\end{align*}
Furthermore, using estimate \eqref{4.D},
\begin{align*}
  J_5 &\le \tau\|\bm{r}_3^{k+1}\|_{-1,2}
  \|\bm{u}_h^{k+1}-\widehat{\bm{u}}_h^{k+1}\|_{1,2} \\
  &\le C(\delta)\tau\|\bm{r}_3^{k+1}\|_{-1,2}^2 
  + \delta\tau\mathcal{D}(0,\bm{u}_h^{k+1}-\widehat{\bm{u}}_h^{k+1})  
  + C\tau\mathcal{E}\big(\pvec{\rho}_h^{k+1},\bm{u}_h^{k+1}|
  \pvec{\widehat\rho}_h^{k+1},\widehat{\bm{u}}_h^{k+1}\big).
\end{align*}
Finally, by the Poincar\'e--Wirtinger inequality for $\vec{\mu}_h^{k+1} - \pvec{\widehat\mu}_h^{k+1}$, we have
\begin{align*}
  J_6 &= -\tau\langle r_1^{k+1}+s_{1}^{k+1}, 
  \vec{\mu}_h^{k+1} - \pvec{\widehat\mu}_h^{k+1} + r_2^{k+1}\rangle \\
  &\leq C(\delta)\tau\|r_1^{k+1}+s_1^{k+1}\|_{-1,2}^2
  + \delta\tau\mathcal{D}(\vec{\mu}_h^{k+1}
  - \pvec{\widehat\mu}_h^{k+1},0) 
  + C\tau\langle \vec{\mu}_h^{k+1} - \pvec{\widehat\mu}_h^{k+1}, \vec{1}\rangle^2 + C\tau\|r_2^{k+1}\|_{1,2}^2.
\end{align*}
It follows from 
\begin{align*}
  \langle \vec{\mu}_h^{k+1} - \pvec{\widehat\mu}_h^{k+1}, 
  \vec{1}\rangle^2 
  &= \sum_{i=1}^N\bigg\langle\frac{\pa f}{\pa\rho_i}
  (\pvec{\rho}_h^{k+1}) - \frac{\pa f}{\pa\rho_i}
  (\pvec{\widehat\rho}_h^{k+1})-r_{2,i}^{k+1},1 
  \bigg\rangle^2 \\
  &\leq C\|\pvec{\rho}_h^{k+1}-\pvec{\widehat\rho}_h^{k+1}\|_{0,2}^2+C\|r_2^{k+1}\|_{0,2}^2
\end{align*}
(here we used \eqref{N2} and \eqref{pertN2} with $\xi_i=1$ and $\langle p_h^{k+1},1\rangle=\langle\widehat p_h^{k+1},1\rangle=0$)
that
\begin{align*}
  J_6 \le C(\delta)\tau\|r_1^{k+1}+s_1^{k+1}\|_{-1,2}^2
  + \delta\tau\mathcal{D}(\vec{\mu}_h^{k+1}
  - \pvec{\widehat\mu}_h^{k+1},0) 
  + C\tau\|\pvec{\rho}_h^{k+1}-\pvec{\widehat\rho}_h^{k+1}\|_{0,2}^2
  + C\tau\|r_2^{k+1}\|_{1,2}^2.
\end{align*}
Summarizing the estimates for $J_1,\ldots,J_6$, we conclude from \eqref{discre11} for sufficiently small $\delta>0$ that
\begin{align}\label{4.aux}
  \mathcal{E}&(\vec{\rho}_h,\bm{u}_h|
  \vec{\widehat\rho}_h,\widehat{\bm{u}}_h)\Big|_{t^k}^{t^{k+1}}  
  + \frac{\tau}{2}\mathcal{D}\big(\vec{\mu}_h^{k+1}
  -\pvec{\widehat\mu}_h^{k+1},
  \bm{u}_h^{k+1}-\widehat{\bm{u}}_h^{k+1}\big) \\
  &\phantom{xx}+ \frac{\tau^2}{2}\langle \rho_h^k,I_h[|\pa_t^\tau
  (\bm{u}_h^{k+1}-\widehat{\bm{u}}_h^{k+1})|^2]\rangle \nonumber \\
  &\le C\tau\mathcal{E}\big(\pvec{\rho}_h^{k+1},\bm{u}_h^{k+1}|
  \pvec{\widehat\rho}_h^{k+1},\widehat{\bm{u}}_h^{k+1}\big) 
  + C\tau\mathcal{E}\big(\pvec{\rho}_h^{k},\bm{u}_h^k|
  \pvec{\widehat\rho}_h^{k},\widehat{\bm{u}}_h^{k}\big) 
  - \tau\langle\bm{s}_3^{k+1},
  \bm{u}_h^{k+1}-\widehat{\bm{u}}_h^{k+1} \rangle \nonumber \\
  &\phantom{xx}+ \tau\bigg(\sum_{i=1}^N\big( \|r_{1,i}^{k+1}+s_{1,i}^{k+1}\|_{-1,2}^2 
  +\|r_{2,i}^{k+1}\|_{1,2}^2\big) 
  + \|\bm{r}_3^{k+1}\|_{-1,2}^2\bigg) + C\tau^3. \nonumber 
\end{align}
The constant $C>0$ depends on the uniform bounds of $\rho_{h,i}^{k}$, $\widehat{\rho}_{h,i}^k$ and the $L^\infty(\Omega)$ bounds for $M_\eps\na\widehat{\mu}^{k+1}$ and $\pa_t^\tau\widehat{\rho}_{h,i}^{k+1}$ for $i=1,\ldots,N$.

In the following, we estimate the term involving $\bm{s}_3^{k+1}$. Using definition \eqref{s3}, we have
\begin{align*}
  & -\langle\bm{s}_3^{k+1},\bm{u}_h^{k+1}-\widehat{\bm{u}}_h^{k+1} \rangle
  = -(J_7 + \cdots + J_{12}), \quad\mbox{where} \\
  & J_7 = \frac{1}{2}\big\langle \pa_t^\tau
  (\rho_h^{k+1}-\widehat\rho_h^{k+1}),I_h[\widehat{\bm{u}}_h^{k+1}
  (\bm{u}_h^{k+1}-\widehat{\bm{u}}_h^{k+1})]\big\rangle, \\
  & J_8 = \big\langle \rho_h^k-\widehat\rho_h^k,
  I_h[\pa_t^\tau\widehat{\bm{u}}_h^{k+1}
  (\bm{u}_h^{k+1}-\widehat{\bm{u}}_h^{k+1})] \big\rangle, \\
  & J_9 =  \bm{b}_{\rm skw}(\rho_h^*\bm{u}_h^*
  - \widehat\rho_h^*\widehat{\bm{u}}_h^*;
  \widehat{\bm{u}}_h^{k+1},\bm{u}_h^{k+1}-\widehat{\bm{u}}_h^{k+1}), \\
  & J_{10} =  \bm{b}^\eps_{\rm skw}(\bm{J}_h^{k+1}-\widehat{\bm{J}}_h^{k+1};
  \widehat{\bm{u}}_h^{k+1},\bm{u}_h^{k+1}-\widehat{\bm{u}}_h^{k+1}), \\
  & J_{11} = \langle \mathbb{S}(\vec{\rho}_{h}^*,
  \na\widehat{\bm{u}}_h^{k+1})
  - \mathbb{S}(\widehat{\vec{\rho}}_{h}^*,
  \na\widehat{\bm{u}}_h^{k+1}),
  \na(\bm{u}_h^{k+1}-\widehat{\bm{u}}_h^{k+1}) \rangle, \\
  & J_{12} = \sum_{i=1}^N\big\langle 
  (\rho_{h,i}^{k}-\widehat\rho_{h,i}^k)
  \na(\widehat\mu_{h,i}^{k+1} - V_i\widehat p_h^{k+1}), 
  \bm{u}_h^{k+1}-\widehat{\bm{u}}_h^{k+1}\big\rangle.
\end{align*}
Some terms can be estimated in a standard way, leading to 
\begin{align*}
   J_8 &\le \|\rho_h^k-\widehat\rho_h^k\|_{0,2}
  \|\pa_t^\tau\widehat{\bm{u}}_h^{k+1}\|_{0,\infty}
  \|I_h[\bm{u}_h^{k+1}-\widehat{\bm{u}}_h^{k+1}]\|_{0,2} \\
  &\le \|\rho_h^k-\widehat\rho_h^k\|_{0,2}
  \|\pa_t^\tau\widehat{\bm{u}}_h^{k+1}\|_{0,\infty}
  (\|\bm{u}_h^{k+1}-\widehat{\bm{u}}_h^{k+1}\|_{0,2}
  + \|\bm{u}_h^{k+1}-\widehat{\bm{u}}_h^{k+1}\|_{1,2}) \\
  &\le C\mathcal{E}\big(\pvec{\rho}_h^{k},\bm{u}_h^k|
  \pvec{\widehat\rho}_h^{k},\widehat{\bm{u}}_h^{k}\big) 
  + C\|\pa_t^\tau\widehat{\bm{u}}_h^{k+1}\|_{0,\infty}^2
  \mathcal{E}\big(\pvec{\rho}_h^{k+1},\bm{u}_h^{k+1}|
  \pvec{\widehat\rho}_h^{k+1},\widehat{\bm{u}}_h^{k+1}\big)\\
  &\phantom{xx}+ \delta\mathcal{D}
  \big(0,\bm{u}_h^{k+1} - \widehat{\bm{u}}_h^{k+1}\big), \\
  J_9 &\le C\big(\|\rho_h^*\|_{0,\infty}
  \|\bm{u}_h^{*} - \widehat{\bm{u}}_h^{*}\|_{0,2} 
  + \|\widehat{\bm{u}}_h^{*}\|_{0,\infty}
  \|\rho_h^*-\widehat\rho_h^{*}\|_{0,2}\big) \\
  &\phantom{xx}\times 
  \big(\|\bm{u}_h^{k+1} - \widehat{\bm{u}}_h^{k+1}\|_{0,2}
  + \|\na(\bm{u}_h^{k+1} - \widehat{\bm{u}}_h^{k+1})\|_{0,2}\big) \\
  &\le\delta\mathcal{D}
  \big(0,\bm{u}_h^{k+1} - \widehat{\bm{u}}_h^{k+1}\big) 
  + C\mathcal{E}\big(\pvec{\rho}_h^{k+1},\bm{u}_h^{k+1}|
  \pvec{\widehat\rho}_h^{k+1},\widehat{\bm{u}}_h^{k+1}\big) 
  + C\mathcal{E}\big(\pvec{\rho}_h^{k},\bm{u}_h^{k}|
  \pvec{\widehat\rho}_h^{k},\widehat{\bm{u}}_h^{k}\big), \\
  J_{11} &\le \delta\mathcal{D}(0,\bm{u}_h^{k+1}
  - \widehat{\bm{u}}_h^{k+1}) 
  + C\|\pvec{\rho}_h^* - \pvec{\widehat{\rho}}_h^*\|_{0,2}^2
  \|\na\widehat{\bm{u}}_h^{k+1}\|_{0,\infty}^2 \\
  &\le \delta\mathcal{D}\big(0,\bm{u}_h^{k+1}
  - \widehat{\bm{u}}_h^{k+1}\big) \\
  &\phantom{xx}+ C\|\na\widehat{\bm{u}}_h^{k+1}\|_{0,\infty}^2
  \big[\mathcal{E}\big(\pvec{\rho}_h^{k},
  \bm{u}_h^k|\pvec{\widehat\rho}_h^{k},\widehat{\bm{u}}_h^{k}\big)
  + \mathcal{E}\big(\pvec{\rho}_h^{k+1},\bm{u}_h^{k+1}|
  \pvec{\widehat\rho}_h^{k+1},\widehat{\bm{u}}_h^{k+1}\big)\big] \\
  J_{12} &= \sum_{i=1}^N\langle (\rho_{h,i}^{k}-\widehat\rho_{h,i}^k)
  \na\widehat\mu_{h,i}^{k+1}, 
  \bm{u}_h^{k+1} - \widehat{\bm{u}}_h^{k+1}\rangle \\
  &\le\delta\mathcal{D}\big(0,\bm{u}_h^{k+1} 
  - \widehat{\bm{u}}_h^{k+1}\big) 
  + C(1+\|\na\widehat\mu_h^{k+1}\|_{0,3}^2)
  \mathcal{E}\big(\pvec{\rho}_h^{k},\bm{u}_h^k|
  \pvec{\widehat\rho}_h^{k},\widehat{\bm{u}}_h^{k}\big) \\
  &\phantom{xx}
  + C\mathcal{E}\big(\pvec{\rho}_h^{k+1},\bm{u}_h^{k+1}|
  \pvec{\widehat\rho}_h^{k+1},\widehat{\bm{u}}_h^{k+1}\big).
\end{align*}

It remains to estimate $J_7$ and $J_{10}$:
\begin{align*}
  J_7 &= \frac{1}{2}\big\langle \pa_t^\tau(\rho_h^{k+1}-\widehat\rho_h^{k+1}),I_h[
  (\bm{u}_h^{k+1}-\widehat{\bm{u}}_h^{k+1})
  \cdot\widehat{\bm{u}}_h^{k+1}]\big\rangle \\
  &= \frac{1}{2}\big\langle\rho_h^k(\bm{u}_h^{k+1}
  - \widehat{\bm{u}}_h^{k+1}),
  \na I_h[(\bm{u}_h^{k+1} - \widehat{\bm{u}}_h^{k+1})
  \cdot\widehat{\bm{u}}_h^{k+1}]\rangle \\
  &\phantom{xx}+ \frac{1}{2}\langle \bm{J}_h^{k+1}-\widehat{\bm{J}}_h^{k+1},
  \na I_h[(\bm{u}_h^{k+1} - \widehat{\bm{u}}_h^{k+1})
  \cdot\widehat{\bm{u}}_h^{k+1}] \rangle \\
  &\phantom{xx}- \frac12\sum_{i=1}^N\langle r_{1,i}^{k+1}+s_{1,i}^{k+1},
  I_h[(\bm{u}_h^{k+1} - \widehat{\bm{u}}_h^{k+1})
  \cdot\widehat{\bm{u}}_h^{k+1}] \rangle \\
  &\le C\|\bm{u}_h^{k+1}-\widehat{\bm{u}}_h^{k+1}\|_{0,2}
  \|\na I_h[(\bm{u}_h^{k+1}-\widehat{\bm{u}}_h^{k+1})
  \cdot\widehat{\bm{u}}_h^{k+1}]\|_{0,2} \\
  &\phantom{xx}+ \frac{1}{2}\langle \bm{J}_h^{k+1}-\widehat{\bm{J}}_h^{k+1},
  \na I_h[(\bm{u}_h^{k+1} - \widehat{\bm{u}}_h^{k+1})
  \cdot\widehat{\bm{u}}_h^{k+1}] \rangle \\
  &\phantom{xx}+ C\|r_1^{k+1}+s_1^{k+1}\|_{-1,2}
  \|I_h[(\bm{u}_h^{k+1} - \widehat{\bm{u}}_h^{k+1})
  \cdot\widehat{\bm{u}}_h^{k+1}]\|_{1,2}.
\end{align*}
The norm can be estimated by means of the $H^1$ stability of the nodal interpolant:
\begin{align*}
  \|&\bm{u}_h^{k+1}-\widehat{\bm{u}}_h^{k+1}\|_{0,2} \|\na I_h[(\bm{u}_h^{k+1}-\widehat{\bm{u}}_h^{k+1})
  \cdot\widehat{\bm{u}}_h^{k+1}]\|_{0,2} \\
  &\leq C\|\bm{u}_h^{k+1}-\widehat{\bm{u}}_h^{k+1}\|_{0,2} 
  \|\na [(\bm{u}_h^{k+1}-\widehat{\bm{u}}_h^{k+1})
  \cdot\widehat{\bm{u}}_h^{k+1}]\|_{0,2} \\
  &\leq \delta\mathcal{D}\big(0,\bm{u}_h^{k+1} 
  - \widehat{\bm{u}}_h^{k+1}\big)  
  + C\big(\|\widehat{\bm{u}}_h^{k+1}\|_{0,\infty}^2
  + \|\widehat{\bm{u}}_h^{k+1}\|_{1,\infty}\big)
  \mathcal{E}\big(\pvec{\rho}_h^{k+1},\bm{u}_h^{k+1}|
  \pvec{\widehat\rho}_h^{k+1},\widehat{\bm{u}}_h^{k+1}\big).
\end{align*}

Recalling $J_{10}$ and addition of suitable zeros yields
\begin{align*}
 J_{10} &= \frac{1}{2}\sum_{i=1}^d 
\bigg\langle \bm{J}_{h,i}^{k+1} - \widehat{\bm{J}}_{h,i}^{k+1}, \frac{\pa\widehat{\bm{u}}_h^{k+1}}{\pa x_i}
\cdot(\bm{u}_h^{k+1}-\widehat{\bm{u}}_h^{k+1}) 
- \frac{\pa}{\pa x_i}(\bm{u}_h^{k+1}-\widehat{\bm{u}}_h^{k+1})
\cdot\widehat{\bm{u}}_h^{k+1}\bigg\rangle \\
& = \sum_{i=1}^d 
\bigg\langle \bm{J}_{h,i}^{k+1} - \widehat{\bm{J}}_{h,i}^{k+1}, \frac{\pa\widehat{\bm{u}}_h^{k+1}}{\pa x_i}
\cdot(\bm{u}_h^{k+1}-\widehat{\bm{u}}_h^{k+1}) 
- \frac{1}{2}\frac{\pa}{\pa x_i}I_h\left[(\bm{u}_h^{k+1}-\widehat{\bm{u}}_h^{k+1})
\cdot\widehat{\bm{u}}_h^{k+1}\right]\bigg\rangle \\
& + \frac{1}{2}\sum_{i=1}^d 
\bigg\langle \bm{J}_{h,i}^{k+1} - \widehat{\bm{J}}_{h,i}^{k+1}, \frac{\pa}{\pa x_i}I_h\left[(\bm{u}_h^{k+1}-\widehat{\bm{u}}_h^{k+1})
\cdot\widehat{\bm{u}}_h^{k+1}\right] - \frac{\pa}{\pa x_i}\left[(\bm{u}_h^{k+1}-\widehat{\bm{u}}_h^{k+1})
\cdot\widehat{\bm{u}}_h^{k+1}\right]\bigg\rangle
\end{align*}
we find after summation that
\begin{align*}
  J_7+J_{10} &\leq \delta\mathcal{D}\big(0,\bm{u}_h^{k+1} 
  - \widehat{\bm{u}}_h^{k+1}\big)  
  + C\big(\|\widehat{\bm{u}}_h^{k+1}\|_{0,\infty}^2
  + \|\widehat{\bm{u}}_h^{k+1}\|_{1,\infty}\big)
  \mathcal{E}\big(\pvec{\rho}_h^{k+1},\bm{u}_h^{k+1}|
  \pvec{\widehat\rho}_h^{k+1},\widehat{\bm{u}}_h^{k+1}\big) \\
  &+ \sum_{i=1}^d\bigg\langle \bm{J}^{k+1}_{h,i}
  -\widehat{\bm{J}}_{h,i}^{k+1}, 
  [\bm{u}_h^{k+1} - \widehat{\bm{u}}_h^{k+1}] 
  \frac{\partial\widehat{\bm{u}}_h^{k+1}}{\partial x_i}\bigg\rangle \\
  &+ \frac{1}{2}\sum_{i=1}^d 
\bigg\langle \bm{J}_{h,i}^{k+1} - \widehat{\bm{J}}_{h,i}^{k+1}, \frac{\pa}{\pa x_i}I_h\left[(\bm{u}_h^{k+1}-\widehat{\bm{u}}_h^{k+1})
\cdot\widehat{\bm{u}}_h^{k+1}\right] - \frac{\pa}{\pa x_i}\left[(\bm{u}_h^{k+1}-\widehat{\bm{u}}_h^{k+1})
\cdot\widehat{\bm{u}}_h^{k+1}\right]\bigg\rangle \\
  &\leq \delta\mathcal{D}\big(0,\bm{u}_h^{k+1} 
  - \widehat{\bm{u}}_h^{k+1}\big)  
  + C\big(\|\widehat{\bm{u}}_h^{k+1}\|_{0,\infty}^2
  + \|\widehat{\bm{u}}_h^{k+1}\|_{1,\infty}\big)
  \mathcal{E}\big(\pvec{\rho}_h^{k+1},\bm{u}_h^{k+1}|
  \pvec{\widehat\rho}_h^{k+1},\widehat{\bm{u}}_h^{k+1}\big) \\
  &\phantom{xx}+ \delta \mathcal{D}
  \big(\vec{\mu}_h^{k+1}-\pvec{\widehat\mu}_h^{k+1},0\big) 
  + C\|\nabla\widehat{\bm{u}}_h^{k+1}\|_{0,\infty}^2
  \mathcal{E}\big(\pvec{\rho}_h^{k+1},\bm{u}_h^{k+1}|
  \pvec{\widehat\rho}_h^{k+1},\widehat{\bm{u}}_h^{k+1}\big) \\
  &+ \frac{1}{2}\sum_{i=1}^d 
\bigg\langle \bm{J}_{h,i}^{k+1} - \widehat{\bm{J}}_{h,i}^{k+1}, \frac{\pa}{\pa x_i}I_h\left[(\bm{u}_h^{k+1}-\widehat{\bm{u}}_h^{k+1})
\cdot\widehat{\bm{u}}_h^{k+1}\right] - \frac{\pa}{\pa x_i}\left[(\bm{u}_h^{k+1}-\widehat{\bm{u}}_h^{k+1})
\cdot\widehat{\bm{u}}_h^{k+1}\right]\bigg\rangle
\end{align*}

It remains to estimate the last line from above, which can be done in an abstract way.
\begin{lemma}
   Let $\bm{a}\in \nabla\mathcal{V}_h$ and $\bm{b},\bm{c}\in\mathcal{X}_h$. Then the following estimate holds with a constant $C$ independent of $\bm{a},\bm{b},\bm{c}$ and $h$:
   \begin{align}
     |\langle \bm{a}, \nabla\left[\bm{b}\cdot\bm{c} - I_h(\bm{b}\cdot\bm{c}) \right] \rangle| \leq C\|\bm{a}\|_{0,2}\|\bm{b}\|_{0,2}\|\bm{c}\|_{1,\infty}. 
   \end{align}
\end{lemma}

\begin{proof}
  We split the integral into element contributions, i.e.
  \begin{equation*}
    \langle \bm{a}, \nabla\left[\bm{b}\cdot\bm{c} - I_h(\bm{b}\cdot\bm{c})\right] \rangle = \sum_{K\in\mathcal{T}_h} \langle \bm{a}, \nabla\left[\bm{b}\cdot\bm{c} - I_h(\bm{b}\cdot\bm{c})\right] \rangle_K  
  \end{equation*}
  For the elementwise contribution we introduce $\tilde{\bm{c}}\in\mathbb{P}^1(K)$ as $\tilde{\bm{c}} = \bm{c} - \bm{c}_K$, with $\bm{c}_K$ is the value at the barycenter of element $K$. For this, we see by a simple Taylor expansion that 
  \begin{align*}
      \tilde{\bm{c}_i} = \nabla \bm{c}_i\vert_K\cdot(x-x_K).
  \end{align*}
  Hence, we can directly estimate
  \begin{align*}
      \|\bm{b}\tilde{\bm{c}}\|_{0,2,K} \leq \|\bm{b}\|_{0,2,K}\|\tilde{\bm{c}}\|_{0,\infty,K} \leq Ch_K\|\bm{b}\|_{0,2,K}\|\nabla\bm{c}\|_{0,\infty,K}.
  \end{align*}
  Furthermore, using the interpolation property, we find
  \begin{align*}
      \bm{b}\cdot\bm{c} - I_h[\bm{b}\cdot\bm{c}] = \bm{b}\cdot\tilde{\bm{c}} - I_h[\bm{b}\cdot\tilde{\bm{c}}]
  \end{align*}
  as the part $\bm{c}_K\bm{b}$ is linear and therefore nodal interpolation is exact.
  Hence, we estimate using $H^1$-stability on $K$ and inverse inequality
  \begin{align*}
     |\langle \bm{a}, \nabla\left[\bm{b}\cdot\bm{c} - I_h(\bm{b}\cdot\bm{c})\right] \rangle_K| &\leq \|\bm{a}\|_{0,2,K}\|\nabla\left[\bm{b}\cdot\tilde{\bm{c}} - I_h[\bm{b}\cdot\tilde{\bm{c}}]\right]\|_{0,2,K} \leq C\|\bm{a}\|_{0,2,K}\|\bm{b}\cdot\tilde{\bm{c}}\|_{1,2,K} \\
     & \leq Ch_K^{-1}\|\bm{a}\|_{0,2,K}\|\bm{b}\cdot\tilde{\bm{c}}\|_{0,2,K}. 
  \end{align*}
  Insertion of our previous estimate and cancellation of the powers of $h_K$ yields
  \begin{align*}
     |\langle \bm{a}, \nabla\left[\bm{b}\cdot\bm{c} - I_h(\bm{b}\cdot\bm{c})\right] \rangle_K| \leq C\|\bm{a}\|_{0,2,K}\|\bm{b}\|_{0,2,K}\|\nabla\bm{c}\|_{0,\infty,K} 
  \end{align*}

  Finally, summation over all elements and application of the Cauchy-Schwarz inequality yields
  \begin{align*}
    |\langle \bm{a}, \nabla\left[\bm{b}\cdot\bm{c} - I_h(\bm{b}\cdot\bm{c})\right] \rangle | \leq C\|\bm{a}\|_{0,2}\|\bm{b}\|_{0,2}\|\bm{c}\|_{1,\infty}.
  \end{align*}
\end{proof}

The result can then be applied to the previous term and yields
\begin{align*}
 &\frac{1}{2}\sum_{i=1}^d 
\bigg\langle \bm{J}_{h,i}^{k+1} - \widehat{\bm{J}}_{h,i}^{k+1}, \frac{\pa}{\pa x_i}I_h\left[(\bm{u}_h^{k+1}-\widehat{\bm{u}}_h^{k+1})
\cdot\widehat{\bm{u}}_h^{k+1}\right] - \frac{\pa}{\pa x_i}\left[(\bm{u}_h^{k+1}-\widehat{\bm{u}}_h^{k+1})
\cdot\widehat{\bm{u}}_h^{k+1}\right]\bigg\rangle \\
&\leq C\|\bm{J}_{h}^{k+1} - \widehat{\bm{J}}_{h}^{k+1}\|_{0,2}\|\widehat{\bm{u}}_h^{k+1}\|_{1,\infty}\|\bm{u}_h^{k+1}-\widehat{\bm{u}}_h^{k+1}\|_{0,2}
\end{align*}
Then using Young's inequality with $\delta>0$, we conclude that
\begin{align*}
  J_7+J_{10} \le \delta\mathcal{D}\big(\vec{\mu}_h^{k+1} 
  - \pvec{\widehat{\mu}}_h^{k+1},
  \bm{u}_h^{k+1} - \widehat{\bm{u}}_h^{k+1}\big) 
  + C(\delta)\mathcal{E}\big(\pvec{\rho}_h^{k+1},\bm{u}_h^{k+1}|
  \pvec{\widehat\rho}_h^{k+1},\widehat{\bm{u}}_h^{k+1}\big).   
\end{align*}
Finally, we estimate the norm of $s_{1,i}^{k+1}$, defined in \eqref{s1}:
\begin{equation*}
  \|s_{1,i}^{k+1}\|_{-1,2}^2
  \le \|(\rho_{h,i}^k-\widehat{\rho}_{h,i}^k)\widehat{\bm{u}}_h^{k+1}
  \|_{0,2}^2
  \le  C\|\widehat{\bm{u}}_h^{k+1}\|_{0,\infty}^2
  \mathcal{E}\big(\pvec{\rho}_h^{k},\bm{u}_h^k|
  \pvec{\widehat\rho}_h^{k},\widehat{\bm{u}}_h^{k}\big).
\end{equation*}

Summarizing the estimates for the terms involving $\bm{s}_3^{k+1}$ and $s_{1,i}^{k+1}$, we conclude from \eqref{4.aux} for sufficiently small $\delta>0$ that
\begin{align}
  \mathcal{E}\big(&\vec{\rho}_h,\bm{u}_h|
  \vec{\widehat\rho}_h,\widehat{\bm{u}}_h\big)\Big|_{t^k}^{t^{k+1}} 
  + \frac{\tau}{2}\mathcal{D}\big(\vec{\mu}_h^{k+1}
  - \pvec{\widehat\mu}_h^{k+1},\bm{u}_h^{k+1}-\widehat{\bm{u}}_h^{k+1}) \label{eq:discrelC} \\
  &\le C\tau\mathcal{E}\big(\pvec{\rho}_h^{k+1},\bm{u}_h^{k+1}|
  \pvec{\widehat\rho}_h^{k+1},\widehat{\bm{u}}_h^{k+1}\big) 
  + C\tau\mathcal{E}\big(\pvec{\rho}_h^{k},\bm{u}_h^k|
  \pvec{\widehat\rho}_h^{k},\widehat{\bm{u}}_h^{k}) \nonumber \\
  &\phantom{xx}+ \tau\bigg(\sum_{i=1}^N
  \big(\|r_{1,i}^{k+1}\|_{-1,2}^2 +\|r_{2,i}^{k+1}\|_{1,2}^2\big)
  + \|\bm{r}_3^{k+1}\|_{-1,2}^2\bigg) + C\tau^3. \nonumber 
\end{align}
Notice that the specific choice of quadrature rules allows the convective terms involving the mass flux to cancel, eliminating a higher-order term. Neither the convective term of the mass flux nor the time derivative can be controlled alone. We observe that when $\eps=0$, the estimate is significantly easier, since the mass flux vanishes in this case.


\subsection{Residual error estimates}

We need to estimate the residuals containing the projection errors:
\begin{align}
  \langle r_{1,i}^{k+1},\psi_i \rangle 
  &= \langle \pa_t^\tau\widehat\rho_{h,i}^{k+1} 
  - \pa_t\rho_i(t^{k+1}),\psi_i\rangle 
  - \langle \widehat\rho_{h,i}^k\widehat{\bm{u}}_h^{k+1}
  - \rho_i(t^{k+1})\bm{u}(t^{k+1}),\na\psi_i \rangle \nonumber \\
  &\phantom{xx}+ \sum_{j=1}^N\langle M_{\eps,ij}(\pvec{\widehat\rho}_h^*)
  \na\widehat\mu_{h,j}^{k+1} - M_{\eps,ij}(\vec{\rho}(t^{k+1}))
  \na\mu_j(t^{k+1}),\na\psi_i \rangle, \label{5.res} \\
  \langle r_{2,i}^{k+1},\xi_i \rangle 
  &= \red{-}\bigg\langle \frac{\pa f}{\pa\rho_i}(\pvec{\widehat\rho}_h^{k+1}) 
  - \frac{\pa f}{\pa\rho_i}(\vec{\rho}(t^{k+1})),\xi_i \bigg\rangle 
  \red{-} V_i\langle \widehat p_h^{k+1}-p(t^{k+1}),\xi_i \rangle, 
  \nonumber \\
  \langle \bm{r}_{3}^{k+1},\bm{v} \rangle 
  &=  \frac12\langle\widehat{\bm{u}}_h^{k+1}
  \pa_t^\tau\widehat\rho_h^{k+1} 
  - \bm{u}(t^{k+1})\pa_t\rho(t^{k+1}),\bm{v} \rangle 
  + \langle \widehat\rho_h^k\pa_t^\tau\widehat{\bm{u}}_h^{k+1} 
  - \rho(t^{k+1})\pa_t \bm{u}(t^{k+1}),\bm{v}\rangle \nonumber \\
    &\phantom{xx}+ \bm{b}_{\rm skw}(\widehat\rho_h^*\widehat{\bm{u}}_h^*;
  \widehat{\bm{u}}_h^{k+1},\bm{v}) 
  - \bm{b}_{\rm skw}(\rho(t^{k+1})\bm{u}(t^{k+1});
  \bm{u}(t^{k+1}),\bm{v}) \nonumber \\
  &\phantom{xx}+ \bm{b}^\eps_{\rm skw}(\widehat{\bm{J}}_h^{k+1};
  \widehat{\bm{u}}_h^{k+1},\bm{v}) 
  - \bm{b}^\eps_{\rm skw}(\bm{J}(t^{k+1});\bm{u}(t^{k+1}),\bm{v}) 
  - \langle\widehat p_h^{k+1} - p(t^{k+1}),\diver\bm{v}\rangle 
  \nonumber \\
  &\phantom{xx} + \langle\mathbb{S}(\pvec{\widehat\rho}_h^*,
  \na\widehat{\bm{u}}_h^{k+1})
  - \mathbb{S}(\pvec{\rho}(t^{k+1}),
  \na\bm{u}(t^{k+1})), \na\bm{v}\rangle \nonumber \\
  &\phantom{xx}+ \sum_{i=1}^N \big\langle\widehat\rho_{h,i}^k
  \na(\widehat\mu_{h,i}^{k+1} - V_i \widehat p_h^{k+1}) 
  - \rho_i(t^{k+1})
  \na\big(\mu_i(t^{k+1})-V_i p(t^{k+1})\big),\bm{v}\big\rangle 
  \nonumber \\
   &\phantom{xx}+ \frac12\langle\pa_t^\tau\widehat\rho_h^{k+1} ,I_h[\widehat{\bm{u}}_h^{k+1}\bm{v}]
   - \widehat{\bm{u}}_h^{k+1}\bm{v} \rangle + \langle \widehat\rho_h^k,I_h[\pa_t^\tau\widehat{\bm{u}}_h^{k+1}\bm{v}] 
  - \pa_t^\tau\widehat{\bm{u}}_h^{k+1}\bm{v}\rangle, \nonumber \\
  \langle r_{4}^{k+1},q \rangle 
  &= \langle \diver(\widehat{\bm{u}}_h^{k+1} 
  - \bm{u}(t^{k+1})),q\rangle \nonumber \\
  &\phantom{xx}+ \sum_{i,j=1}^N \big\langle 
  M_{\eps,ij}(\pvec{\widehat\rho}_h^*)
  \na\widehat\mu_{h,j}^{k+1} - M_{\eps,ij}(\vec{\rho}(t^{k+1}))
  \na\mu_j(t^{k+1}),V_i\na q \big\rangle. \nonumber 
\end{align}

\begin{lemma}\label{lem.res}
The residuals can be estimated for $i=1,\ldots,N$ as
\begin{align*}
  \|r_{1,i}\|_{-1,2}^2 &\le C(h^4\tau^{-2} + \tau^2 + h^2), \\
  \|r_{2,i}\|_{1,2}^2 &\leq C(\tau^2 + h^2), \\
  \|\bm{r}_{3}\|_{-1,2}^2 &\leq C(h^4\tau^{-2} + \tau^2 + h^2).
\end{align*}
\end{lemma}

\begin{proof}
{\em Step 1: estimate of $r_{1,i}^{k+1}$.} We write $\|r_{1,i}\|_{-1,2}^2 \le K_{11}+K_{12}+K_{13}$, where
\begin{align*}
  K_{11} &= \|\pa_t^\tau\widehat\rho_{h,i}^{k+1} 
  - \pa_t\rho_i(t^{k+1})\|_{-1,2}^2, \\
  K_{12} &= \| \widehat\rho_{h,i}^k\widehat{\bm{u}}_h^{k+1}
  - \rho_i(t^{k+1})\bm{u}(t^{k+1})\|_{0,2}^2, \\
  K_{13} &= \sum_{j=1}^N \|M_{\eps,ij}(\pvec{\widehat\rho}_h^*)
  \na\widehat\mu_{h,j}^{k+1} - M_{\eps,ij}(\vec{\rho}(t^{k+1}))
  \na\mu_j(t^{k+1})\|_{0,2}^2.
\end{align*}
The projection error \eqref{eq:projtime} implies that 
\begin{align*}
  K_{11} \le C(h^4\tau^{-2}+\tau^2).
\end{align*}
Next, by the projection error \eqref{eq:projH1} with $m=0$, the time shift error \eqref{eq:projshift}, and the stability bounds \eqref{eq:projstab},
\begin{align*}
  K_{12} &\le \|\widehat{\bm{u}}_h^{k+1}\|^2_{0,\infty}
  \|\widehat\rho_{h,i}^k - \rho_i(t^{k+1})\|_{0,2}^2 
  + \|\rho_i(t^{k+1})\|_{0,\infty}^2
  \|\widehat{\bm{u}}_h^{k+1} - \bm{u}(t^{k+1})\|_{0,2}^2 \\
  &\le C(h^4 + \tau^2).
\end{align*}
Finally, using the projection errors \eqref{eq:projH1} and \eqref{eq:projH1mup} with $m=1$,
\begin{align*}
  K_{13} &\le C\sum_{j=1}^N \big(
  \|\na(\widehat\mu_{h,j}^{k+1} - \mu_j(t^{k+1}))\|_{0,2}^2 
  + \|M_{\eps,ij}(\pvec{\widehat\rho}_h^*) 
  - M_{\eps,ij}(\vec{\rho}(t^{k+1}))\|_{0,2}^2\big) \\
  &\le C(h^2 + \tau^2).  
\end{align*}
This proves the bound for $r_{1,i}^{k+1}$.

{\em Step 2: Estimate for $r_{2,i}^{k+1}$.} We estimate $\|r_{2,i}^{k+1}\|_{1,2}^2 \le K_{21}+K_{22}$, where
\begin{align*}
  K_{21} = \bigg\|\frac{\pa f}{\pa\rho_i}
  (\pvec{\widehat\rho}_h^{k+1}) 
  - \frac{\pa f}{\pa\rho_i}(\pvec{\rho}(t^{k+1}))\bigg\|_{1,2}^2, \quad
  K_{22} = V_i^2\|\widehat p_h^{k+1} - p(t^{k+1})\|_{1,2}^2.
\end{align*}
The projection error \eqref{eq:projH1mup} for the pressure with $m=1$ immediately implies that $K_{22}\le Ch^2$, while the projection error \eqref{eq:projH1} with $m=1$ and the stability estimate \eqref{eq:projstab} (together with the Lipschitz continuity of the second derivative $f''$) yield
\begin{align*}
  K_{21} \le C\|\pvec{\widehat\rho}_h^{k+1}
  - \pvec{\rho}(t^{k+1})\|_{0,2}^2 
  + C\|\na(\pvec{\widehat\rho}_h^{k+1} - \vec{\rho}(t^{k+1}))\|_{0,2}^2
  \le Ch^2.
\end{align*} 

{\em Step 3: Estimate for $\bm{r}_3^{k+1}$.} The estimation of $\bm{r}_3^{k+1}$ is more involved. First, we write $\|\bm{r}_3^{k+1}\|_{-1,2}^2\le K_{31}+\cdots+K_{40}$, where
\begin{align*}
  K_{31} &= \frac{1}{2}\|\widehat{\bm{u}}_h^{k+1}\pa_t^\tau\widehat\rho_h^{k+1} 
  - \bm{u}(t^{k+1})\pa_t\rho(t^{k+1})\|_{-1,2}^2, \\ 
  K_{32} &= \|\widehat\rho_h^k\pa_t^\tau\widehat{\bm{u}}_h^{k+1} 
  - \rho(t^{k+1})\pa_t\bm{u}(t^{k+1})\|_{-1,2}^2, \\
  K_{33} &= \|\bm{b}_{\rm skw}(\widehat\rho_h^*\widehat{\bm{u}}_h^*;
  \widehat{\bm{u}}_h^{k+1},\cdot) 
  - \bm{b}_{\rm skw}(\rho(t^{k+1})\bm{u}(t^{k+1});
  \bm{u}(t^{k+1}),\cdot)\|_{-1,2}^2, \\
  K_{34} &= \|\bm{b}^\eps_{\rm skw}(\widehat{\bm{J}}_h^{k+1};
  \widehat{\bm{u}}_h^{k+1},\cdot) 
  - \bm{b}^\eps_{\rm skw}(\bm{J}(t^{k+1});
  \bm{u}(t^{k+1}),\cdot)\|_{-1,2}^2, \\
  K_{35} & = \|\widehat{p}^{k+1}_h-p(t^{k+1})\|_{0,2}^2, \\
  K_{36} &= \|\mathbb{S}(\pvec{\widehat\rho}_h^*,
  \na\widehat{\bm{u}}_h^{k+1}) 
  - \mathbb{S}(\vec{\rho}(t^{k+1}),
  \na\bm{u}(t^{k+1}))\|_{0,2}^2, \\ 
  K_{37} &= \sum_{i=1}^N \|\widehat\rho_{h,i}^k
  \na\widehat\mu_{h,i}^{k+1} 
  - \rho_i(t^{k+1})\na\mu_i(t^{k+1})\|_{-1,2}^2, \\
  K_{38} & = \frac12\|\pa_t^\tau\widehat\rho_h^{k+1} \big(I_h[\widehat{\bm{u}}_h^{k+1}\cdot]
  - \widehat{\bm{u}}_h^{k+1}\cdot\big) \|_{-1,2}^2 , \\
  K_{39} & =  \| \widehat\rho_h^k\big(
  I_h[\pa_t^\tau\widehat{\bm{u}}_h^{k+1}\cdot] 
  - \pa_t^\tau\widehat{\bm{u}}_h^{k+1}\cdot\big)\|_{-1,2}^2, \\
  K_{40} &= \sum_{i=1}^N\big\|V_i\big(\widehat\rho_{h,i}^k
  \na\widehat p_h^{k+1}-\rho_i(t^{k+1})\na p(t^{k+1})
  \big)\big\|_{-1,2}^2.
\end{align*}

To estimate the products in the $H^{-1}(\Omega)$ norm, we need an auxiliary result. Let $a\in H^{-1}(\Omega)$, $b\in W^{1,3}(\Omega)\cap L^\infty(\Omega)$, and $v\in H^1(\Omega)$. Then $bv\in H^1(\Omega)$ (in up to three space dimensions) and
\begin{align*}
  \|ab\|_{-1,2} &= \sup_{\|v\|_{1,2}=1}\langle ab,v\rangle 
  = \sup_{\|v\|_{1,2}=1}\langle a,bv\rangle
  \le \sup_{\|v\|_{1,2}=1}\|a\|_{-1,2}\|bv\|_{1,2} \\
  &\le \sup_{\|v\|_{1,2}=1}\|a\|_{-1,2}
  \big(\|b\|_{0,\infty}\|v\|_{1,2}
  + \|\na b\|_{0,3}\|v\|_{0,6}\big) \\
  &\le C\|a\|_{-1,2}\big(\|b\|_{0,\infty} + \|\na b\|_{0,3}\big).
\end{align*}
Hence, using the projection errors \eqref{eq:projH1} and \eqref{eq:projtime},
\begin{align*}
  K_{31} &= \frac12\|\widehat{\bm{u}}_h^{k+1}(\pa_t^\tau\widehat\rho_h^{k+1} 
  - \pa_t\rho(t^{k+1})) + \pa_t\rho(t^{k+1})(\widehat{\bm{u}}_h^{k+1}
  - \bm{u}(t^{k+1}))\|_{-1,2}^2 \\
  &\le C\big(\|\na\widehat{\bm{u}}_h^{k+1}\|_{0,3}^2
  + \|\widehat{\bm{u}}_h^{k+1}\|_{0,\infty}^2\big)
  \|\pa_t^\tau\widehat\rho_h^{k+1} 
  - \pa_t\rho(t^{k+1})\|_{-1,2}^2 \\
  &\phantom{xx}+ C\|\pa_t\rho(t^{k+1})\|_{0,\infty}^2
  \|\widehat{\bm{u}}_h^{k+1} - \bm{u}(t^{k+1})\|_{0,2}^2 \\
  &\le C(h^4\tau^{-2} + \tau^2 + h^4).   
\end{align*}
The second term $K_{32}$ is estimated in a similar way, leading to
\begin{align*}
  K_{32} \le C(h^4\tau^{-2} + \tau^2 + h^4).
\end{align*}
The terms $K_{33}$, $K_{35}$, $K_{36}$, and $K_{37}$ are estimated in a similar way, and we end up with
\begin{align*}
  K_{33} &\le \|\widehat\rho^*_h\widehat{\bm{u}}^*_h 
  - \rho(t^{k+1})\bm{u}(t^{k+1})\|_{0,2}^2
  (\|\widehat{\bm{u}}_h^{k+1}\|_{0,3}^2 
  + \|\na\widehat{\bm{u}}_h^{k+1}\|_{0,3}^2) \\
  &\phantom{xx}+ \| \rho(t^{k+1})\bm{u}(t^{k+1})\|_{0,3}^2
  \big(\|\widehat{\bm{u}}_h^{k+1} - \bm{u}(t^{k+1})\|_{0,6}^2
  + \|\na(\widehat{\bm{u}}_h^{k+1} - \bm{u}(t^{k+1}))\|_{0,2}^2\big) \\
  &\le C(\tau^2 + h^2 + h^4), \\
  K_{35} &\le Ch^4, \\
  K_{36} &\le \|\mathbb{S}(\pvec{\widehat\rho}_h^*,
  \na\widehat{\bm{u}}_h^{k+1})
  - \mathbb{S}(\pvec{\widehat\rho}_h^*,
  \na\bm{u}(t^{k+1}))\|_{0,2}^2 \\
  &\phantom{xx}+ \|\mathbb{S}(\pvec{\widehat\rho}_h^*,
  \na\bm{u}(t^{k+1}))
  - \mathbb{S}(\vec{\rho}(t^{k+1}),
  \na\bm{u}(t^{k+1}))\|_{0,2}^2 \\
  & \leq C(h^2 + \tau^2), \\
  K_{37} &\le \|\widehat\rho_{h,i}^k\|_{1,3}^2
  \|\na(\widehat\mu_{h,i}^{k+1} - \mu_i(t^{k+1}))\|_{-1,2}^2 
  + C\|\na\mu_i(t^{k+1})\|_{0,3}^2
  \|\widehat\rho_{h,i}^k - \rho_i(t^{k+1})\|_{0,2}^2 \\ 
  & \leq C(h^4+\tau^2).
\end{align*}

It remains to consider the term $K_{34}$, written as
\begin{align*}
  & K_{34} = \sup_{\|v\|_{1,2}=1}|Z|^2, \quad
  \mbox{where} \\
  & Z = \bm{b}^\eps_{\rm skw}(\widehat{\bm{J}}_h^{k+1};
  \widehat{\bm{u}}_h^{k+1},\bm{v}) 
  - \bm{b}^\eps_{\rm skw}(\bm{J}(t^{k+1});\bm{u}(t^{k+1}),\bm{v}).
\end{align*}
We add and subtract the term $\bm{b}^\eps_{\rm skw}(\bm{J}(t^{k+1});\widehat{\bm{u}}_h^{k+1},\bm{v})$, which leads after a computation to $Z=Z_1+Z_2$, where
\begin{align*}
  Z_1 &=  \bm{b}^\eps_{\rm skw}(\widehat{\bm{J}}_h^{k+1}-\bm{J}(t^{k+1});
  \widehat{\bm{u}}_h^{k+1},\bm{v}), \\
  Z_2 &= \bm{b}^\eps_{\rm skw}(\bm{J}(t^{k+1});
  \widehat{\bm{u}}_h^{k+1}-\bm{u}(t^{k+1}),\bm{v}).
\end{align*}
 Let $\|\bm{v}\|_{1,2}=1$. The terms $Z_1$ and $Z_2$ are estimated in a standard way: 
\begin{align*}
  Z_1 &\le C\|\widehat{\bm{J}}_{h}^{k+1} - \bm{J}(t^{k+1})\|_{0,2}
  \|\widehat{\bm{u}}_h^{k+1}\|_{1,\infty}\|\bm{v}\|_{1,2}
  \le C\eps \sum_{i=1}^N \|\na(\widehat\mu_{h,i}^{k+1} 
  - \mu_i(t^{k+1}))\|_{0,2}\|\bm{v}\|_{1,2} \\
  &\le C\eps h\|\bm{v}\|_{1,2} \le Ch, \\
  Z_2 &\le  C\|\bm{J}(t^{k+1})\|_{0,\infty}
  \|\widehat{\bm{u}}_h^{k+1} 
  - \bm{u}(t^{k+1})\|_{1,2}\|\bm{v}\|_{1,2} 
  \le C\eps \|\widehat{\bm{u}}_h^{k+1}
  - \bm{u}(t^{k+1})\|_{1,2}\|\bm{v}\|_{1,2}
  \le Ch.
\end{align*}
This shows that $Z\le Ch$ and consequently $K_{34}\le Ch^2$.

Finally we consider $K_{38}$ and $K_{39}$, which contain the quadrature errors. We estimate
\begin{align*}
  \sqrt{2K_{38}} &= \sup_{\|v\|_{1,2}=1}\big\langle
  \pa_t^\tau\widehat\rho_h^{k+1}, I_h[\widehat{\bm{u}}_h^{k+1}\bm{v}]
  - \widehat{\bm{u}}_h^{k+1}\bm{v} \big\rangle \\
  &\le \|\pa_t^\tau\widehat\rho_h^{k+1}\|_{0,\infty}
  \sup_{\|v\|_{1,2}=1} \|I_h[\widehat{\bm{u}}_h^{k+1}\bm{v}]
  - \widehat{\bm{u}}_h^{k+1}\bm{v}\|_{0,1} \\
  &\le Ch^2\|\pa_t^\tau\widehat\rho_h^{k+1}\|_{0,\infty}
  \sup_{\|v\|_{1,2}=1} \|\nabla \widehat{\bm{u}}_h^{k+1}\|_{0,2}
  \|\nabla \bm{v}\|_{0,2} \leq Ch^2, \\
  \sqrt{K_{39}} & \le Ch^2\|\widehat\rho_h^{k}\|_{0,\infty} \|\nabla\pa_t^\tau \widehat{\bm{u}}_h^{k+1}\|_{0,2} \\
  &\le Ch^2\tau^{-1}\|\widehat\rho_h^{k}\|_{0,\infty} 
  \big(\|\nabla \widehat{\bm{u}}_h^{k+1}\|_{0,2} +\|\nabla \widehat{\bm{u}}_h^{k}\|_{0,2}\big) \leq Ch^2\tau^{-1}, \\
 \sqrt{K_{40}} 
  &= \|\na(\widehat p_h^{k+1}-p(t^{k+1}))\|_{-1,2}
  \le \|\widehat p_h^{k+1}-p(t^{k+1})\|_{0,2} \le Ch^2.
\end{align*}
Summarizing the estimates for $K_{31},\ldots,K_{40}$, we conclude the bound for $\bm{r}_3^{k+1}$.
\end{proof}


\subsection{Gronwall argument}

With the residual error estimates from Lemma \ref{lem.res}, we can write \eqref{eq:discrelC} under the condition $h=c\tau$ as
\begin{align*}
  \mathcal{E}\big(&\pvec{\rho}_h,\bm{u}_h|
  \pvec{\widehat\rho}_h,\widehat{\bm{u}}_h\big)\Big|_{t^k}^{t^{k+1}}  
  + \frac{\tau}{2}\mathcal{D}
  \big(\vec{\mu}_h^{k+1} - \pvec{\widehat\mu}_h^{k+1},
  \bm{u}_h^{k+1} - \widehat{\bm{u}}_h^{k+1}\big) \\
  &\le C_1\tau\mathcal{E}\big(\pvec{\rho}_h^{k+1},\bm{u}_h^{k+1}|
  \pvec{\widehat\rho}_h^{k+1},\widehat{\bm{u}}_h^{k+1}\big) 
  + C_2\tau\mathcal{E}\big(\pvec{\rho}_h^{k},\bm{u}_h^k|
  \pvec{\widehat\rho}_h^{k},\widehat{\bm{u}}_h^{k}\big) + C_3\tau^3 
  + C_4\tau(h^2 + \tau^2).
\end{align*}
The constants $C_1,\ldots,C_4$ depend only on the a priori bounds of the discrete solution, the uniform bounds of the regular solution, and the lower bounds for $\rho_{h,i}^k$ for $i=1,\ldots,N$. We sum over all time 
steps and choose $\tau\le \tau_0:=1/(2C_1)$ to find, by the discrete Gronwall lemma, that 
\begin{align}\label{eq:discerror}
  \max_{k=1,\ldots,n}\mathcal{E}\big(\pvec{\rho}_h^k,\bm{u}_h^k|
  \pvec{\widehat\rho}_h^k,\widehat{\bm{u}}_h^k\big) 
  &+ \frac{\tau}{2}\sum_{k=1}^{n}\mathcal{D}
  \big(\vec{\mu}_h^{k} - \pvec{\widehat\mu}_h^{k},
  \bm{u}_h^{k} - \widehat{\bm{u}}_h^{k}\big) \\
  &\le C(T)(h^2+\tau^2) + \mathcal{E}\big(\pvec{\rho}_h^0,\bm{u}_h^0|
  \pvec{\widehat\rho}_h^0,\widehat{\bm{u}}_h^0). \nonumber 
\end{align}
Since the initial data match by assumption, the last term vanishes. Using the lower bounds for the viscosities and $M_{\eps}$, we obtain an error estimate for the evolution error.

The relative energy dissipation only gives information on $\na(\vec{\mu}_h^{k+1}-\vec{\widehat\mu}_h^{k+1})$. However, repeating the computations as for $J_6$ in Section \ref{sec.drei} and using the Poincar\'e--Wirtinger inequality, we can extend this bound to a full $H^1(\Omega)$ control.

\subsection{Pressure error estimates}

The relative energy method does not provide a discrete pressure error estimate. A common approach to address this is to leverage the inf-sup stability to derive such estimates. The advantage of the quasi-incompressible system is that pressure estimates can be obtained from the chemical potential by choosing the test function $\xi_i=p_h^{k+1}-\widehat p_h^{k+1}$ in \eqref{N2} and \eqref{pertN2}:
\begin{align*}
  V_i^2\|p_h^{k+1}-\widehat p_h^{k+1}\|_{0,2}^2 
  \leq C\bigg\|\frac{\pa f}{\pa\rho_i}(\pvec{\rho}_h^{k+1})
  - \frac{\pa f}{\pa\rho_i}(\pvec{\widehat\rho}_h^{k+1})\bigg\|_{0,2}^2 
  + \|\mu_{h,i}^{k+1}-\widehat\mu_{h,i}^{k+1}\|_{0,2}^2 {}+ \|r_{2,i}^{k+1}\|_{0,2}^2.   
\end{align*}
Summation over all time steps and the use of error estimate \eqref{eq:discerror} yields
\begin{align*}
  \sum_{k=1}^{n}\tau \|p_h^{k}-\widehat p_h^{k}\|_{0,2}^2 
  \leq C(h^2+\tau^2).  
\end{align*}
We combine this estimate with the standard projection errors to infer the error estimate of Theorem \ref{thm.error}. 
In the case $\eps=0$, the above argument has to be adjusted by considering the difference $\mu_i-\mu_j$, where $V_i\neq V_j$, since we do not obtain an estimate for all $\mu_i$.


\section{Discussion of the case $\eps=0$}\label{sec.eps}

We outline the difficulties that arise in the case $\eps=0$, when the diffusion matrix is not positive definite. While most terms can be treated as before, the estimation of $r_{1,i}^{k+1}+s_{1,i}^{k+1}$ is more subtle. We need to estimate these terms in some discrete norm with negative Sobolev index; otherwise, we cannot expect any convergence order. We start as follows:
\begin{align*}
  \langle\pvec{r}_1^{k+1}+\pvec{s}_1^{k+1},
  \vec{\mu}_h^{k+1}-\pvec{\widehat\mu}_h^{k+1} \rangle 
  &= \big\langle \pvec{r}_1^{k+1}+\pvec{s}_1^{k+1},
  (\mathrm{I}-\mathbb{P}_{1})
  \big[\vec{\mu}_h^{k+1}-\pvec{\widehat\mu}_h^{k+1}\big] \big\rangle \\ 
  &\phantom{xx}+ \langle\pvec{r}_1^{k+1}+\pvec{s}_1^{k+1},
  \mathbb{P}_{1}[\vec{\mu}_h^{k+1}-\pvec{\widehat\mu}_h^{k+1}] \rangle,
\end{align*}
where $\mathbb{P}_{1}$ is the projection into the subspace generated by the one-vector. The first term can be treated as before, since $M\approx \mathrm{I}-\mathbb{P}_{1}$. Hence, only the projection into the kernel of $M$ needs to be considered:
\begin{align*}
  &N\big\langle \pvec{r}_1^{k+1}+\pvec{s}_1^{k+1},\mathbb{P}_{1}
  \big[\vec{\mu}_h^{k+1}-\pvec{\widehat\mu}_h^{k+1}\big] \big\rangle \\
  &= \langle \mathrm{1}\cdot(\pvec{r}_1^{k+1}+\pvec{s}_1^{k+1}), 
  1\cdot(\vec{\mu}_h^{k+1}-\pvec{\widehat\mu}_h^{k+1})\rangle \\
  &= \bigg\langle  \sum_{i=1}^N(r_{1,i}^{k+1} + s_{1,i}^{k+1}),
  \sum_{i=1}^N\bigg(\frac{\pa f}{\pa\rho_i}(\pvec{\rho}_h^{k+1}) 
  - \frac{\pa f}{\pa\rho_i}(\pvec{\widehat\rho}_h^{k+1}) 
  + V_i(p_h^{k+1}-\widehat p_h^{k+1})\bigg)\bigg\rangle.
\end{align*}
The difficulty with this expression lies in the fact that estimates in negatively indexed Sobolev norms appear out of reach, as we lack gradient bounds for both the partial densities and the pressure.

Another route would be to consider $M$ as a kind of vanishing viscosity approximation, i.e.\ $M_h=M + h^\beta\mathrm{I}$ with $\beta>0$. It turns out that the estimate of $\langle s_{1,i}^{k+1}, \mu_{h,i}^{k+1} - \widehat\mu_{h,i}^{k+1}\rangle$ is problematic. We may reformulate this term in the continuous case, for a compressible variant, as in \cite[Appendix E.1]{Dru23}:
\begin{align*}
  \bigg\langle\diver(\widehat\rho_i\widehat{\bm{u}}), 
  \frac{\pa f}{\pa\rho_i}(\pvec{\rho}_h^{k+1}) 
  - \frac{\pa f}{\pa\rho_i}(\pvec{\widehat\rho}_h^{k+1}) 
  - \sum_{j=1}^N\frac{\pa^2 f}{\pa\rho_i\pa\rho_j}(\pvec{\widehat\rho}_h^{k+1})
  (\rho_j-\widehat\rho_j)\bigg\rangle.
\end{align*}
For our discretization, a discrete analogue of this formulation cannot be achieved for two reasons. First, we treat $\rho_i$ in this term explicitly, while the contributions from $\mu_i$ are implicit. Second, even if the first problem can be dealt with, we would only have terms of the form $\diver((\rho_i^k-\widehat{\rho_i}^k)\widehat{\bm{u}}_h^{k+1})$ and $\pi_h^0[(\pa f/\pa\rho_i)(\pvec{\rho}_h^{k+1}) - (\pa f/\pa\rho_i)(\pvec{\widehat\rho}_h^{k+1})]$. Hence, we would have to deal with the projection error. This topic is left for future research.

\section{Numerical experiments}\label{sec.exp}

In this section, we present some numerical simulations in two space dimensions, illustrating the convergence of the scheme. The scheme is implemented in NGSolve \cite{Sch14}, using Netgen for the generation of an unstructured mesh of triangles. The resulting nonlinear system is solved by Newton's method with a tolerance of $10^{-9}$ (measured in the $L^2(\Omega)$ norm of the Newton updates), and the linear systems are computed via direct LU decomposition. In all tests, we fixed the computational domain $\Omega=(0,1)^2$ with periodic boundary conditions. This means that we identify $\Omega$ with the two-dimensional torus. 

\subsection{Convergence tests}

We consider $N=2$ components, the final time $T=0.1$, and the stress tensor $\mathbb{S}(\na\bm{u})=\nu(\na\bm{u}^T+\na\bm{u}) + \lambda(\diver\bm{u})\mathbb{I}$ with $\nu=10^{-3}$ and $\lambda=0$. The internal energy density is chosen as
$f(\vec\rho)=\sum_{i=1}^2\rho_i\log(\rho_i/\rho)$, and the mobility matrix is given by $M_{ij}(\pvec{\rho})=(\rho_i+\eps)\delta_{ij} - \rho_i\rho_j/\rho$. The initial data equal
\begin{align*}
  & \rho_1^0(x,y) = 1+0.8\sin(4\pi x)\sin(2\pi y), \quad
  \rho_2^0(x,y) = V_2^{-1}(1-V_1\rho_1^0(x,y)), \\
  & \bm{u}^0(x,y) = \big(-\sin(\pi x)^2\sin(2\pi y),
  \sin(2\pi x)\sin(\pi y)^2\big)^T,
\end{align*}
where $V_1=0.3$, $V_2=0.7$.
Since no exact solution is available, we compute the error by using the reference solution at the finest space resolution. Given the mesh sizes $h_k=2^{-k}$ for $k=2,\ldots,7$ and $h_{\rm ref}=2^{-8}$, we define the error quantities 
\begin{align*}
  & \mbox{err}(\vec\rho) = \max_{\ell=1,\ldots,n} \sqrt{
  \sum_{i=1}^2\|\rho_{h_k,i}^\ell-\rho_{{\rm ref},i}^\ell
  \|_{L^2(\Omega)}^2}, \quad
  \mbox{err}(\vec\mu) = \sqrt{\tau\sum_{\ell=1}^{n}\sum_{i=1}^2
  \|\mu_{h_k,i}^\ell-\mu_{{\rm ref},i}^\ell\|_{L^2(\Omega)}^2}, \\
  & \mbox{err}(\bm{u}) = \max_{\ell=1,\ldots,n}
  \sqrt{\|\bm{u}_{h_k}^\ell-\bm{u}_{\rm ref}^\ell\|_{L^2(\Omega)}^2}, \qquad\quad
  \mbox{err}(p) = \sqrt{\tau\sum_{\ell=1}^{n}\|p_{h_k}^\ell-p_{\rm ref}^\ell
  \|_{L^2(\Omega)}^2}.
\end{align*}
We have chosen the time step size $\tau=10^{-3}$ such that $n=T/\tau=100$. Tables \ref{table_10-2}--\ref{table_0} present the errors and experimental convergence order for $\eps=10^{-2}$, $\eps=10^{-6}$, $\eps=10^{-10}$, and $\eps=0$, respectively. Notice that the tables show only the spatial error, since $\tau$ is fixed and the reference solution is computed with the same $\tau$; the coupling $h=c\tau$ of Theorem \ref{thm.error} is not tested. Moreover, $\mbox{err}(\vec\mu)$ is measured in $L^2(0,T;L^2(\Omega))$, while the theorem gives a bound in $L^2(0,T;H^1(\Omega))$.

\begin{figure}[htbp!]
\centering
\footnotesize
\begin{tabular}{c@{}c}
    \multicolumn{2}{c}{\includegraphics[trim={13cm 33.5cm 9cm 0cm},clip,scale=0.25]{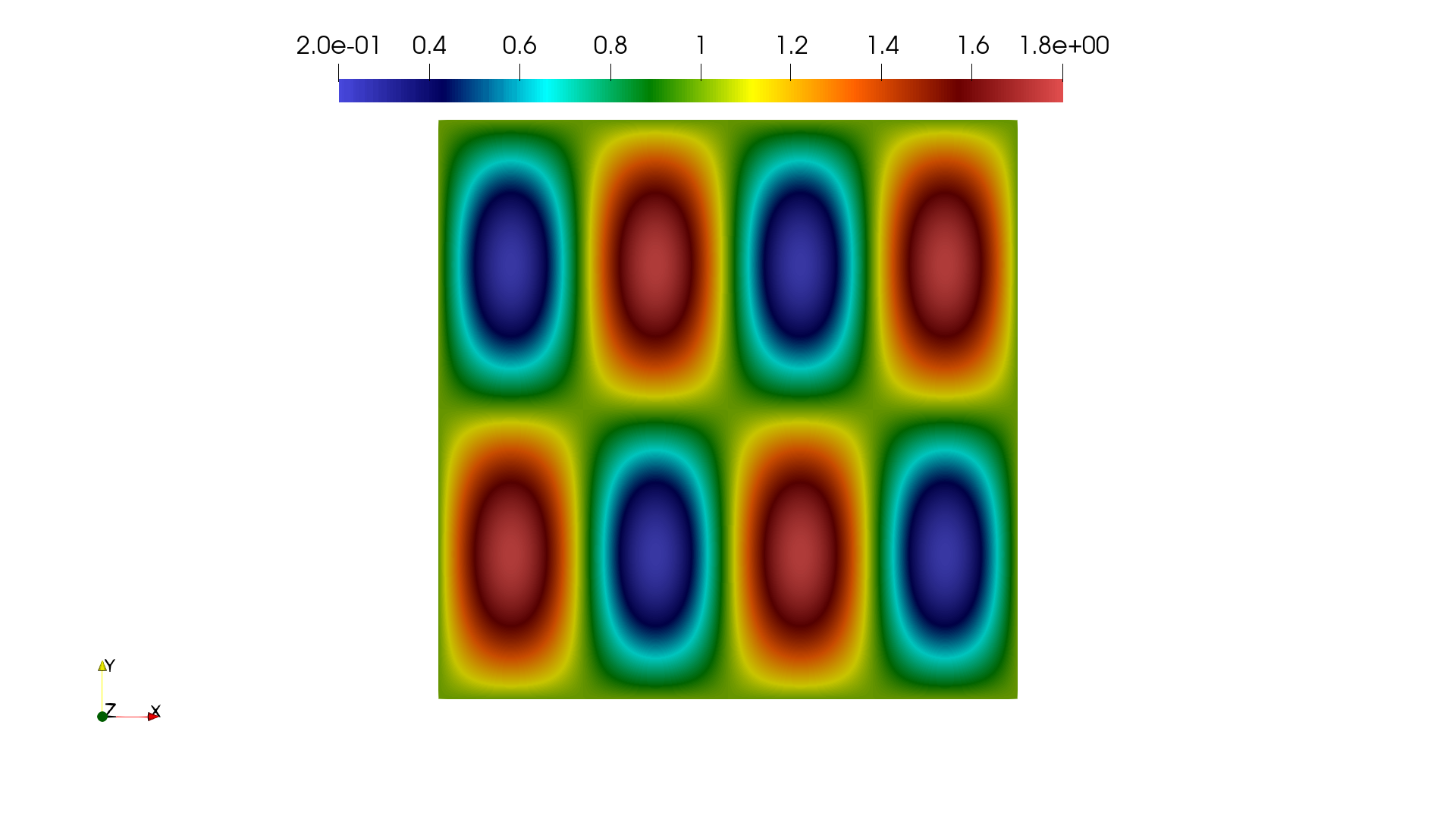}} \\
    \includegraphics[trim={20cm 5cm 20cm 5cm},clip,scale=0.15]{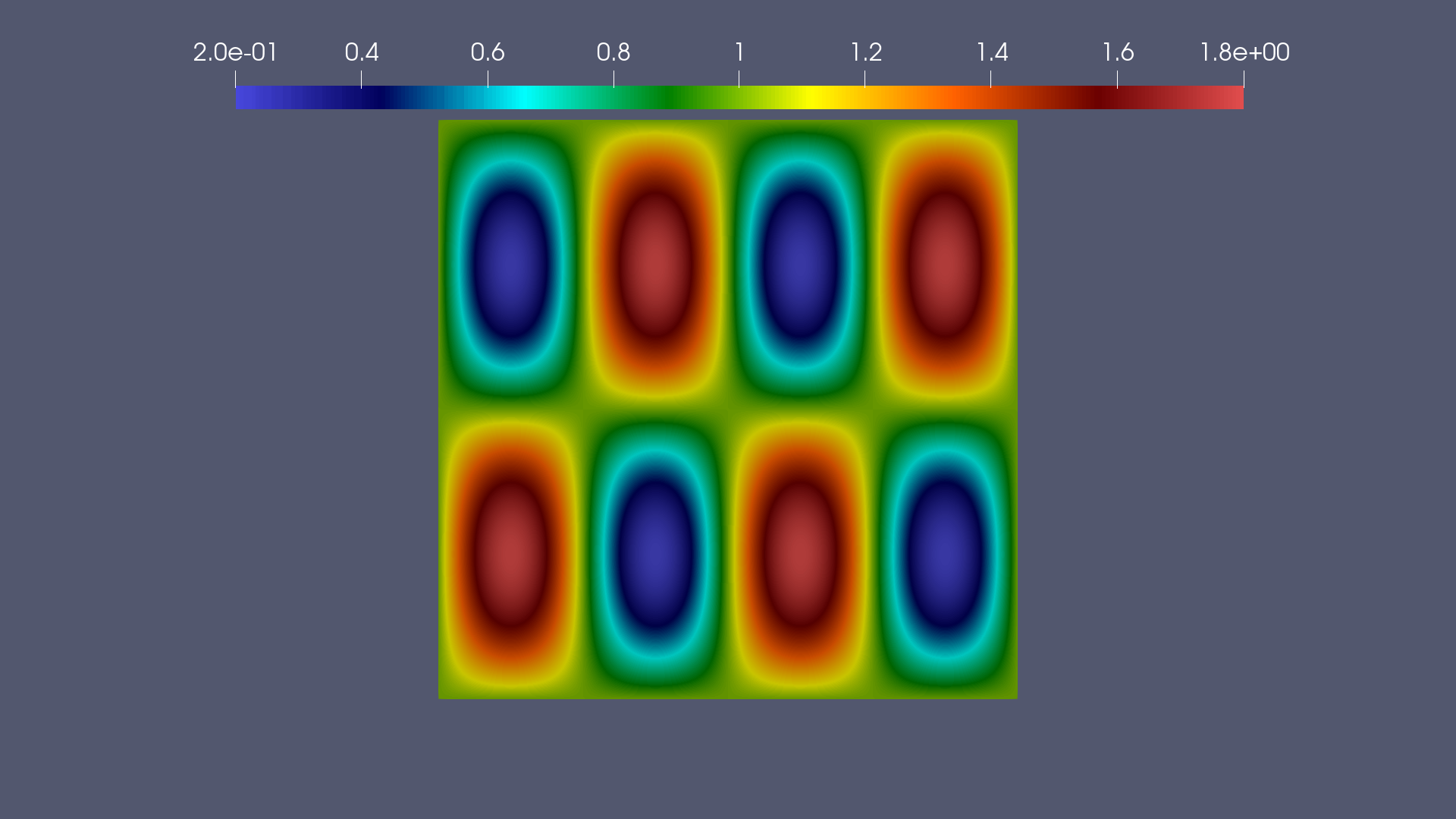} 
    &
    \includegraphics[trim={20cm 5cm 20cm 5cm},clip,scale=0.15]{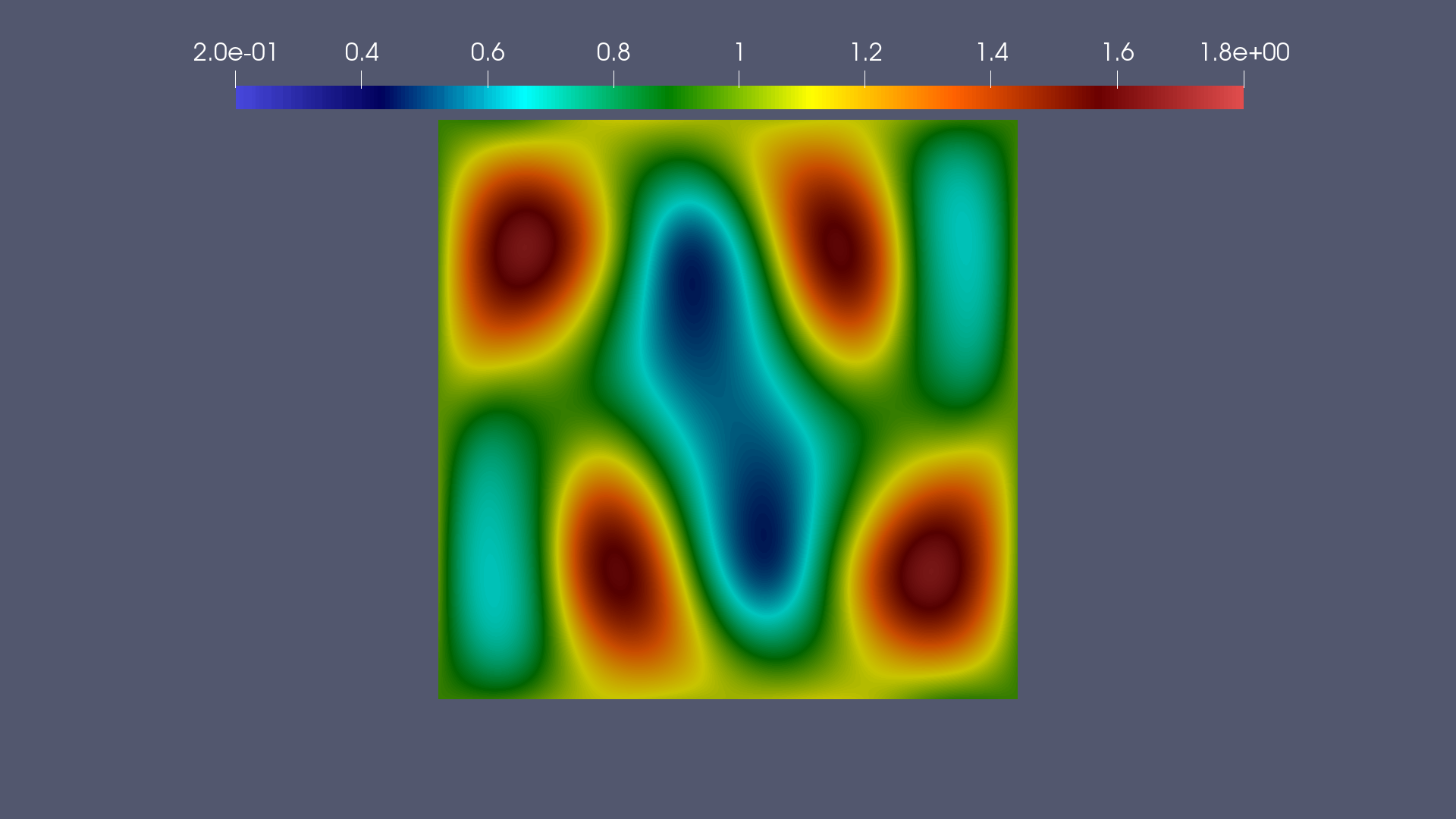}  
     \\
     \multicolumn{2}{c}{\includegraphics[trim={13cm 33.5cm 9cm 0cm},clip,scale=0.25]{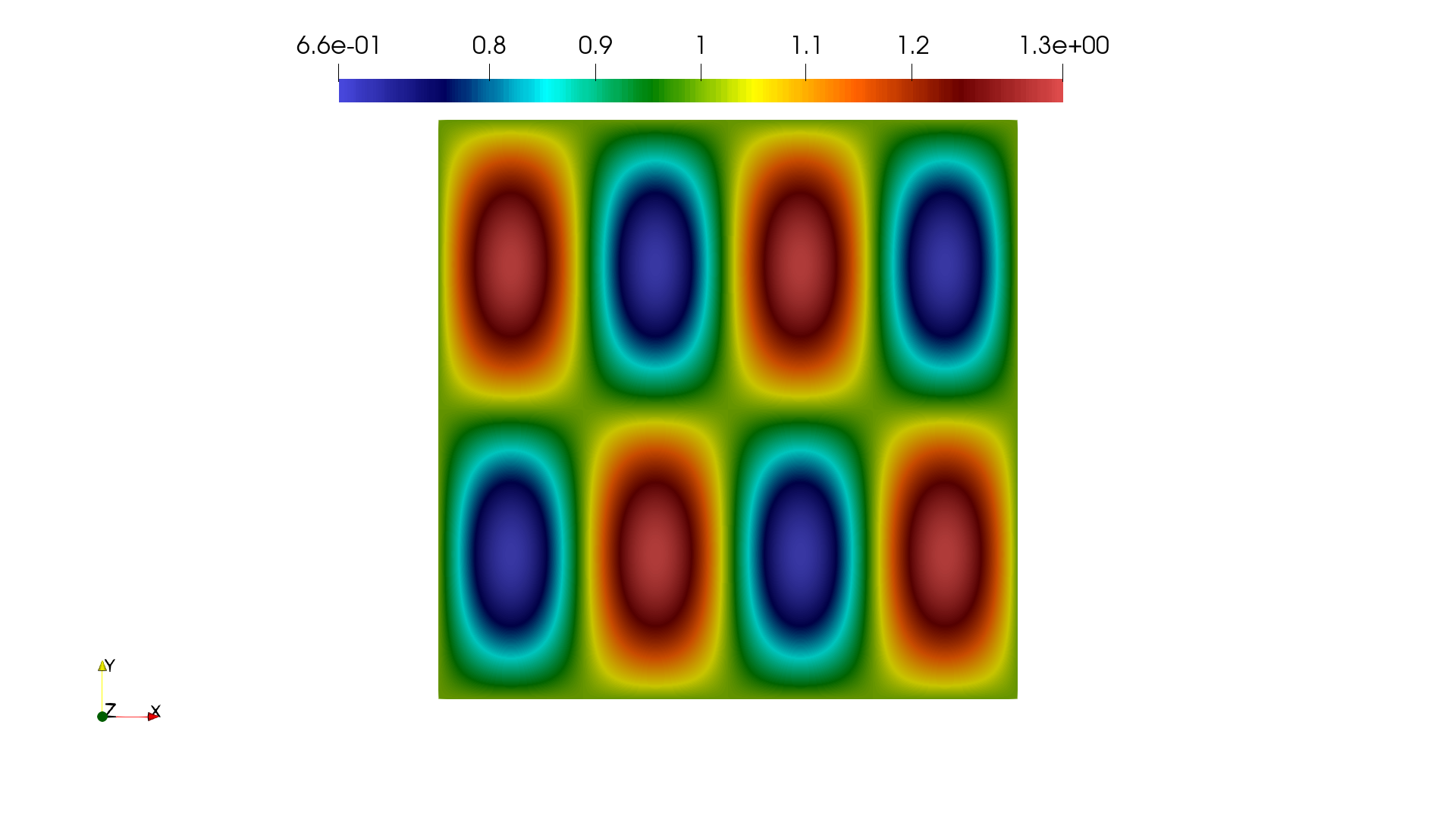}} \\
    \includegraphics[trim={20cm 5cm 20cm 5cm},clip,scale=0.15]{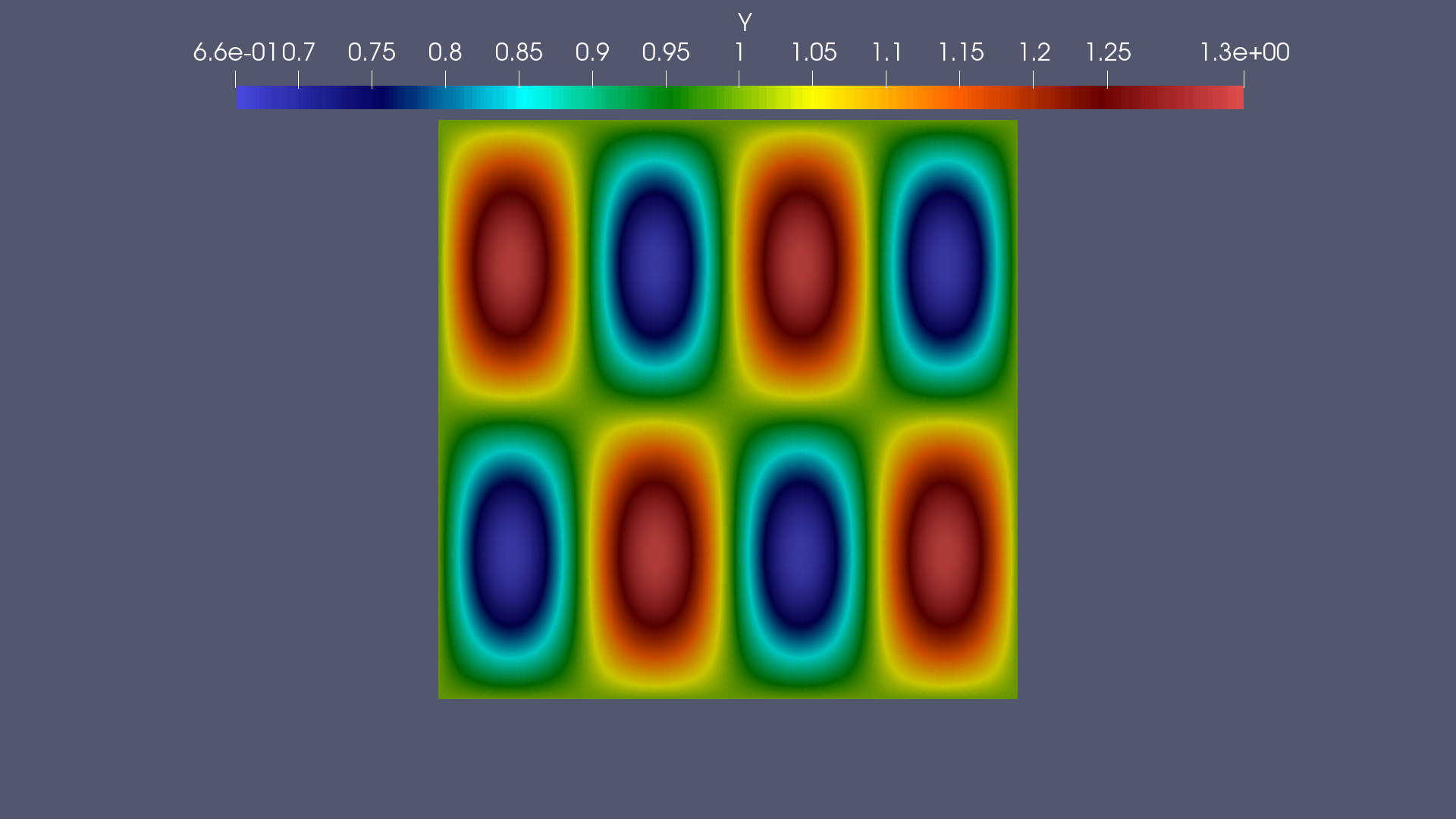}  
    &
    \includegraphics[trim={20cm 5cm 20cm 5cm},clip,scale=0.15]{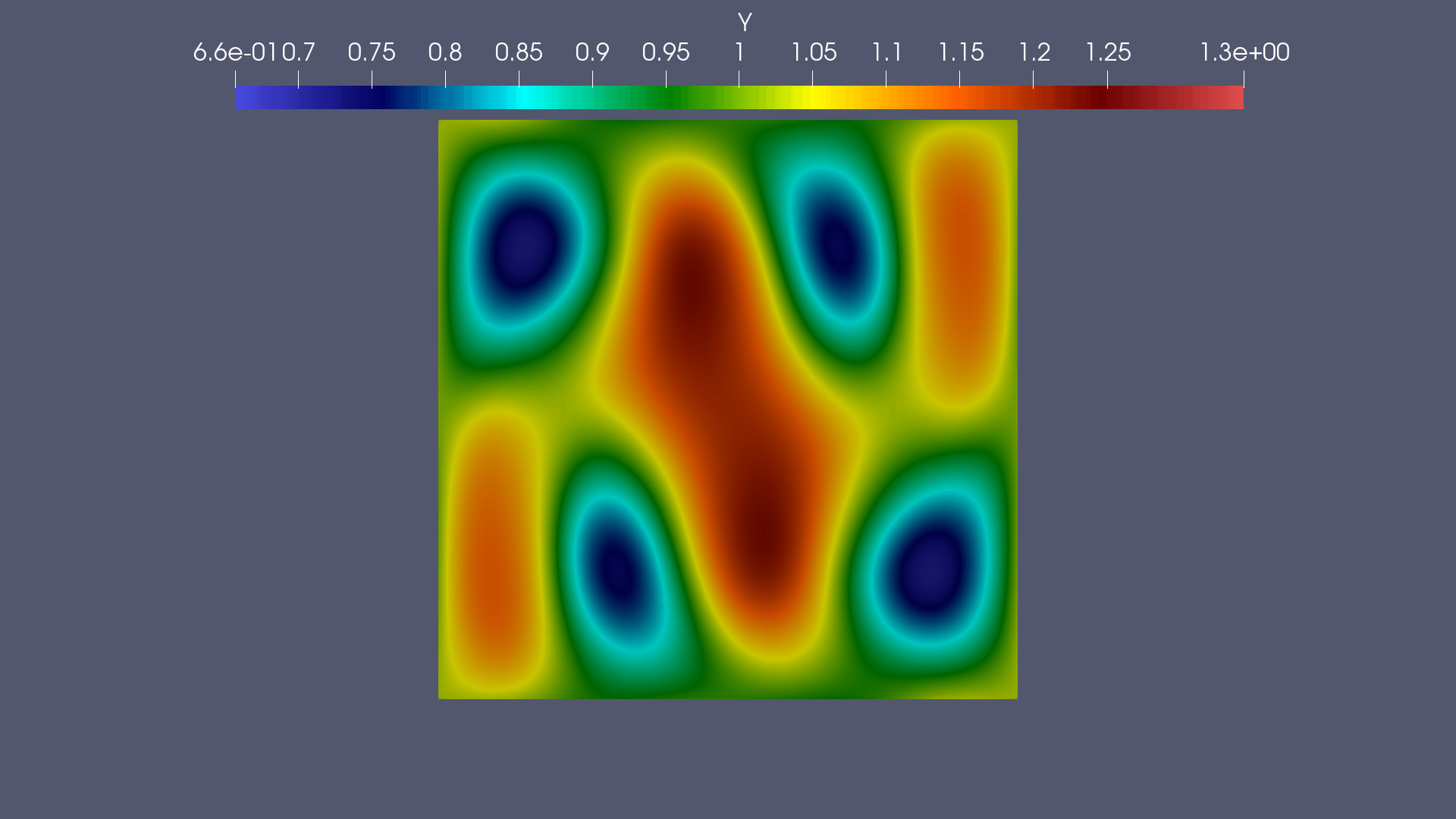} 
     \\
     \multicolumn{2}{c}{\includegraphics[trim={14cm 33.5cm 9cm 0cm},clip,scale=0.25]{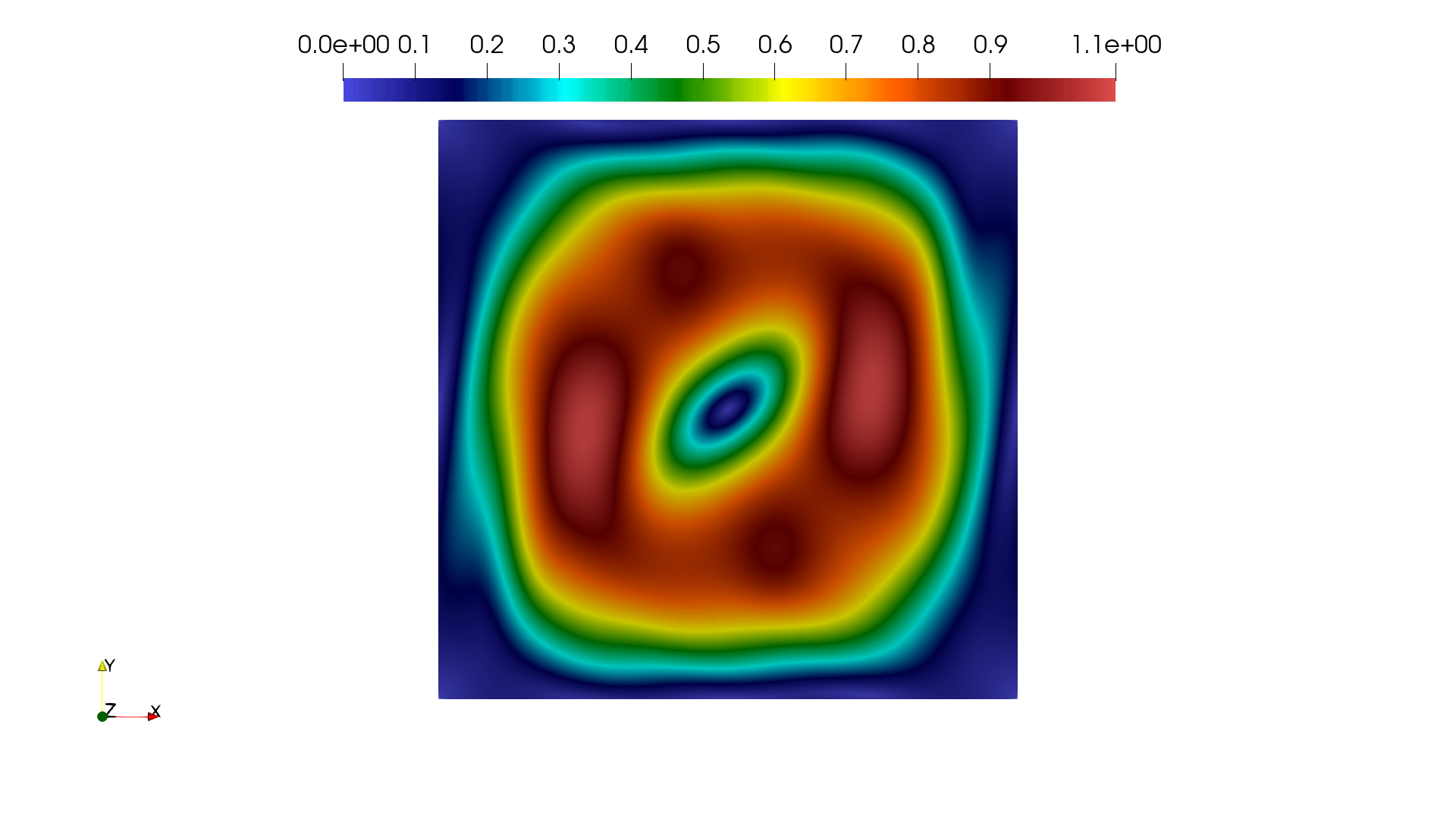}} \\
     \includegraphics[trim={20cm 5cm 20cm 5cm},clip,scale=0.15]{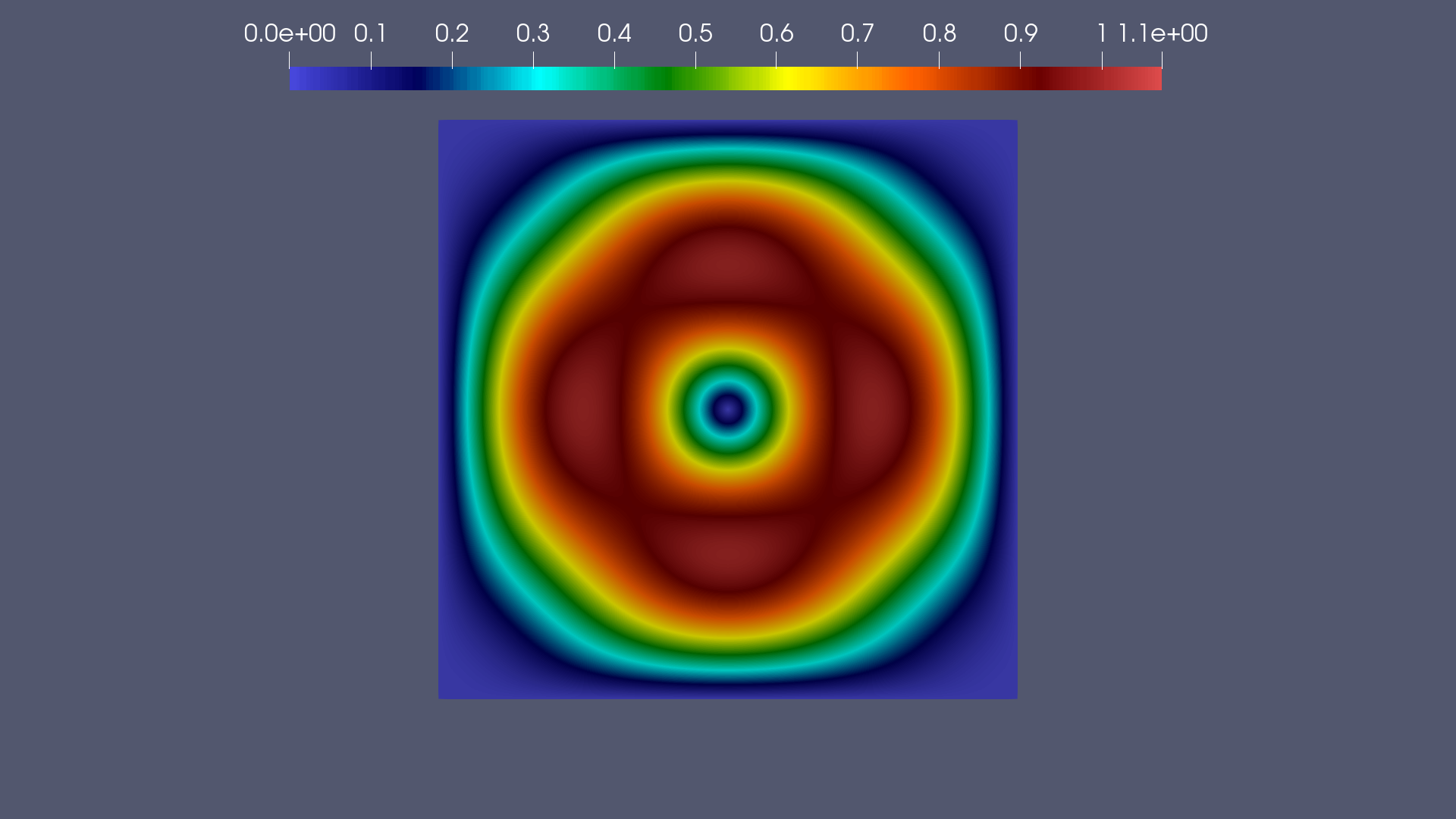}  
    &
    \includegraphics[trim={20cm 5cm 20cm 5cm},clip,scale=0.15]{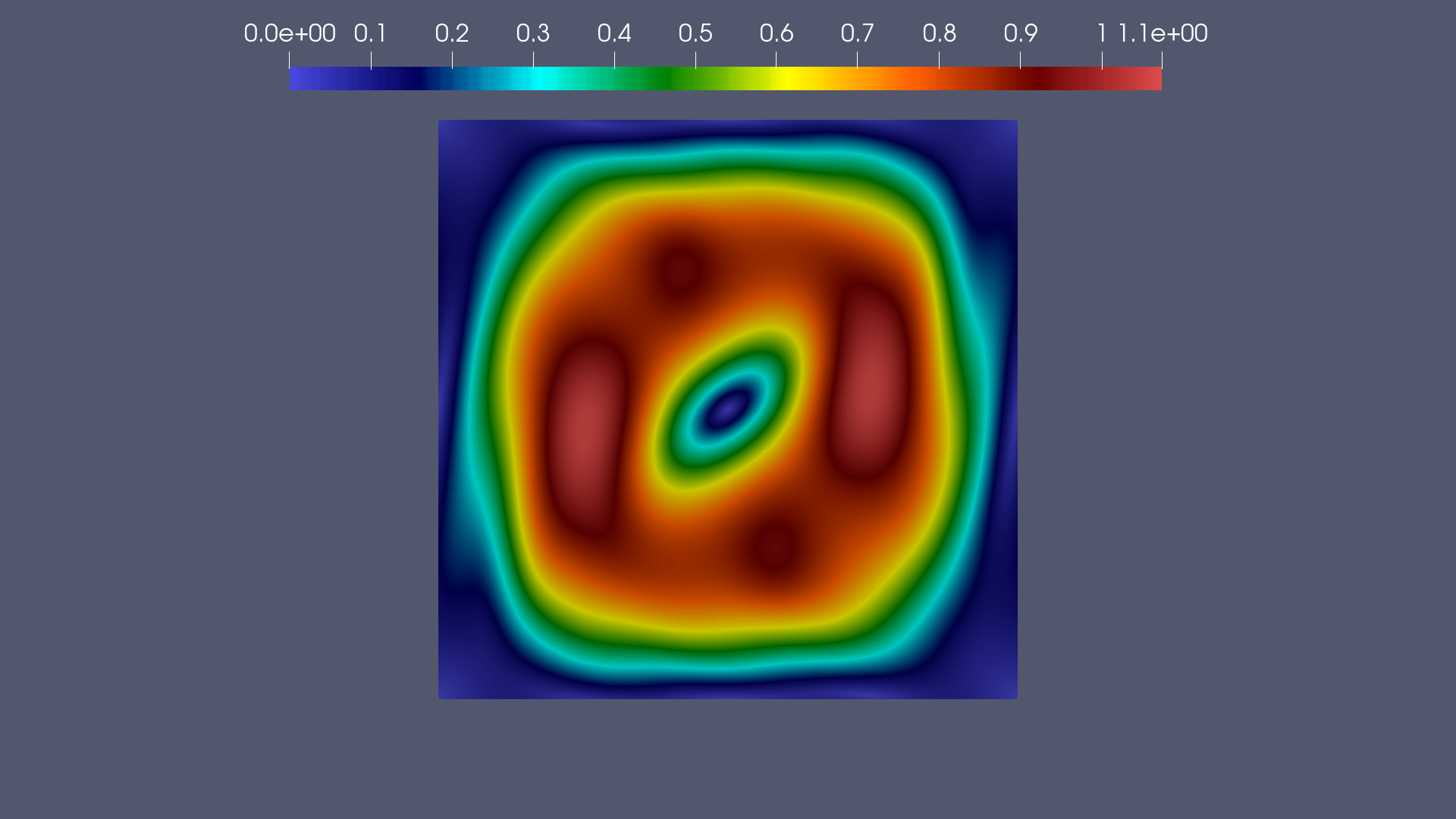} 
     \\
\end{tabular}
    \caption{Snapshots of partial mass densities $\rho_1$ (top), $\rho_2$ (middle) and magnitude of the velocity field $|\bm{u}|$ (bottom) at initial time $t=0$ and final time $t=0.1$.\label{fig:plots}}
\end{figure}

A visualization of the initial and final configurations for $\varepsilon=10^{-10}$ is shown in Figure \ref{fig:plots}.

\begin{table}[htbp!]
	\centering
	\small
	\caption{$L^2(\Omega)$ errors and experimental order of convergence for $\varepsilon=10^{-2}$.}
	\begin{tabular}{|c||c|c|c|c|c|c|c|c|}
		\hline
		$ k $ & $\text{err}(\vec\rho)$ & eoc & $\text{err}(\vec\mu)$ & eoc & $\text{err}(\bm{u})$ & eoc & $\text{err}(p)$ & eoc \\
		\hline
		2 & $3.42\cdot 10^{-1}$ &  -- & $1.20\cdot 10^{-1}$ & -- & $1.82\cdot 10^{-1}$ &  -- & $1.66\cdot 10^{-1}$ &  --\\
		3 & $9.38\cdot 10^{-2}$ & $\phm1.87$ & $5.30\cdot 10^{-2}$ & $\phm1.18$ & $7.75\cdot 10^{-2}$ & $\phm1.24$ & $7.45\cdot 10^{-2}$ & $\phm1.16$ \\
		4 & $2.55\cdot 10^{-2}$ & $\phm1.88$ & $1.49\cdot 10^{-2}$ & $\phm1.83$ & $2.54\cdot 10^{-2}$ & $\phm1.61$ & $2.06\cdot 10^{-2}$ & $\phm1.85$ \\
		5 & $6.48\cdot 10^{-3}$ & $\phm1.98$ & $3.78\cdot 10^{-3}$ & $\phm1.98$ & $7.24\cdot 10^{-3}$ & $\phm1.81$ & $5.20\cdot 10^{-3}$ & $\phm1.99$ \\
		6 & $1.63\cdot 10^{-3}$ & $\phm1.99$ & $9.48\cdot 10^{-4}$ & $\phm2.00$ & $1.86\cdot 10^{-3}$ & $\phm1.96$ & $1.30\cdot 10^{-3}$ & $\phm2.00$ \\
		7 & $4.07\cdot 10^{-4}$ & $\phm2.00$ & $2.37\cdot 10^{-4}$ & $\phm2.00$ & $4.68\cdot 10^{-4}$ & $\phm1.99$ & $3.26\cdot 10^{-4}$ & $\phm2.00$ \\
		\hline
	\end{tabular}
	\label{table_10-2}
\end{table}

\begin{table}[htbp!]
	\centering
	\small
	\caption{$L^2(\Omega)$ errors and experimental order of convergence for $\varepsilon=10^{-6}$.}
	\begin{tabular}{|c||c|c|c|c|c|c|c|c|}
		\hline
		$ k $ & $\text{err}(\vec\rho)$ & eoc & $\text{err}(\vec\mu)$ & eoc & $\text{err}(\bm{u})$ & eoc & $\text{err}(p)$ & eoc \\
		\hline
		2 & $3.42\cdot 10^{-1}$ &  -- & $1.75\cdot 10^{-1}$ &  -- & $2.14\cdot 10^{-1}$ &  -- & $2.42\cdot 10^{-1}$ &  -- \\
		3 & $1.83\cdot 10^{-1}$ & $\phm0.90$ & $1.02\cdot 10^{-1}$ & $\phm0.79$ & $1.09\cdot 10^{-1}$ & $\phm0.97$ & $1.44\cdot 10^{-1}$ & $\phm0.75$ \\
		4 & $4.97\cdot 10^{-2}$ & $\phm1.88$ & $3.01\cdot 10^{-2}$ & $\phm1.76$ & $4.71\cdot 10^{-2}$ & $\phm1.22$ & $4.15\cdot 10^{-2}$ & $\phm1.79$ \\
		5 & $1.32\cdot 10^{-2}$ & $\phm1.91$ & $8.36\cdot 10^{-3}$ & $\phm1.85$ & $1.30\cdot 10^{-2}$ & $\phm1.86$ & $1.13\cdot 10^{-2}$ & $\phm1.88$ \\
		6 & $3.28\cdot 10^{-3}$ & $\phm2.01$ & $2.13\cdot 10^{-3}$ & $\phm1.98$ & $3.18\cdot 10^{-3}$ & $\phm2.03$ & $2.84\cdot 10^{-3}$ & $\phm1.99$ \\
		7 & $8.14\cdot 10^{-4}$ & $\phm2.01$ & $5.26\cdot 10^{-4}$ & $\phm2.02$ & $7.83\cdot 10^{-4}$ & $\phm2.02$ & $7.03\cdot 10^{-4}$ & $\phm2.01$ \\
		\hline
	\end{tabular}
	\label{table_10-6}
\end{table}

\begin{table}[htbp!]
	\centering
	\small
	\caption{$L^2(\Omega)$ errors and experimental order of convergence for $\varepsilon=10^{-10}$.}
	\begin{tabular}{|c||c|c|c|c|c|c|c|c|}
		\hline
		$ k $ & $\text{err}(\vec\rho)$ & eoc & $\text{err}(\vec\mu)$ & eoc & $\text{err}(\bm{u})$ & eoc & $\text{err}(p)$ & eoc \\
		\hline
		2 & $3.42\cdot 10^{-1}$ &  -- & $1.75\cdot 10^{-1}$ &  -- & $2.14\cdot 10^{-1}$ &  -- & $2.42\cdot 10^{-1}$ &  -- \\
		3 & $1.83\cdot 10^{-1}$ & $\phm0.90$ & $1.02\cdot 10^{-1}$ & $\phm0.79$ & $1.09\cdot 10^{-1}$ & $\phm0.97$ & $1.44\cdot 10^{-1}$ & $\phm0.75$ \\
		4 & $4.97\cdot 10^{-2}$ & $\phm1.88$ & $3.01\cdot 10^{-2}$ & $\phm1.76$ & $4.71\cdot 10^{-2}$ & $\phm1.22$ & $4.15\cdot 10^{-2}$ & $\phm1.79$ \\
		5 & $1.32\cdot 10^{-2}$ & $\phm1.91$ & $8.37\cdot 10^{-3}$ & $\phm1.85$ & $1.30\cdot 10^{-2}$ & $\phm1.86$ & $1.13\cdot 10^{-2}$ & $\phm1.88$ \\
		6 & $3.28\cdot 10^{-3}$ & $\phm2.01$ & $2.13\cdot 10^{-3}$ & $\phm1.98$ & $3.18\cdot 10^{-3}$ & $\phm2.03$ & $2.84\cdot 10^{-3}$ & $\phm1.99$ \\
		7 & $8.15\cdot 10^{-4}$ & $\phm2.01$ & $5.26\cdot 10^{-4}$ & $\phm2.02$ & $7.83\cdot 10^{-4}$ & $\phm2.02$ & $7.03\cdot 10^{-4}$ & $\phm2.02$ \\
		\hline
	\end{tabular}
	\label{table_10-10}
\end{table}

\begin{table}[htbp!]
	\centering
	\small
	\caption{$L^2(\Omega)$ errors and experimental order of convergence for $\varepsilon=0$.}
	\begin{tabular}{|c||c|c|c|c|c|c|c|c|}
		\hline
		$ k $ & $\text{err}(\vec\rho)$ & eoc & $\text{err}(\vec\mu)$ & eoc & $\text{err}(\bm{u})$ & eoc & $\text{err}(p)$ & eoc \\
		\hline
		2 & $3.42\cdot 10^{-1}$ &  -- & $1.75\cdot 10^{-1}$ &  -- & $2.14\cdot 10^{-1}$ &  -- & $2.42\cdot 10^{-1}$ &  -- \\
		3 & $1.83\cdot 10^{-1}$ & $\phm0.90$ & $1.02\cdot 10^{-1}$ & $\phm0.79$ & $1.09\cdot 10^{-1}$ & $\phm0.97$ & $1.44\cdot 10^{-1}$ & $\phm0.75$ \\
		4 & $4.97\cdot 10^{-2}$ & $\phm1.88$ & $3.01\cdot 10^{-2}$ & $\phm1.76$ & $4.71\cdot 10^{-2}$ & $\phm1.22$ & $4.15\cdot 10^{-2}$ & $\phm1.79$ \\
		5 & $1.32\cdot 10^{-2}$ & $\phm1.91$ & $8.37\cdot 10^{-3}$ & $\phm1.85$ & $1.30\cdot 10^{-2}$ & $\phm1.86$ & $1.13\cdot 10^{-2}$ & $\phm1.88$ \\
		6 & $3.28\cdot 10^{-3}$ & $\phm2.01$ & $2.13\cdot 10^{-3}$ & $\phm1.98$ & $3.18\cdot 10^{-3}$ & $\phm2.03$ & $2.84\cdot 10^{-3}$ & $\phm1.99$ \\
		7 & $8.15\cdot 10^{-4}$ & $\phm2.01$ & $5.26\cdot 10^{-4}$ & $\phm2.02$ & $7.83\cdot 10^{-4}$ & $\phm2.02$ & $7.03\cdot 10^{-4}$ & $\phm2.02$ \\
		\hline
	\end{tabular}
	\label{table_0}
\end{table}

We observe second-order convergence for $(\vec\rho,\vec\mu,\bm{u},p)$ for all cases. In all tests, the partial and total masses and the pointwise quasi-incompressibility constraint are conserved up to errors of order $10^{-14}$. Tables \ref{table_10-2}--\ref{table_0} show that the error is insensitive to values of $\varepsilon$. For $\varepsilon \in\{10^{-6},10^{-10},0\}$, the differences in the errors are smaller than $10^{-10}$. 


\end{document}